\documentclass[11pt,reqno]{amsart}
\usepackage[a4paper,margin=2.5cm]{geometry}
\usepackage{amsmath,amssymb,amsthm,mathtools,microtype}
\usepackage[numbers,sort&compress]{natbib}
\usepackage[hidelinks]{hyperref}
\numberwithin{equation}{section}
\newtheorem{theorem}{Theorem}[section]
\newtheorem{lemma}[theorem]{Lemma}
\newtheorem{proposition}[theorem]{Proposition}
\newtheorem{corollary}[theorem]{Corollary}
\theoremstyle{definition}
\newtheorem{definition}{Definition}[section]
\newcommand{\R}{\mathbb R}
\newcommand{\cG}{\mathcal G}
\newcommand{\cZ}{\mathcal Z}

\allowdisplaybreaks
\title[Segregated normalized solutions with two mass constraints]
{Segregated normalized solutions for a nonlinear Schr\"odinger system with two mass constraints}

\author{Qing Guo}
\address{Department of Mathematics, Minzu University of China,
Beijing 100081, China}
\email{guoqing0117@163.com}
\thanks{}
\author{Wenjing Li}
\address{Department of Mathematics, Minzu University of China,
Beijing 100081, China}
\email{wenjingli@muc.edu.cn}
\author{Shiyang Song}
\address{Department of Mathematics, Minzu University of China,
Beijing 100081, China}
\email{shiyangsong@muc.edu.cn}

\date{}

\begin{document}

\begin{abstract}
We study a two-component nonlinear Schr\"odinger system with two prescribed positive $L^2$ masses and a fixed competitive coupling $\beta<0$, for $N=1,2,3$ and $4\le p<2^*$. We first consider families whose peak heights diverge and whose components converge uniformly to zero away from disjoint finite peak sets. For such families, we prove that both Lagrange multipliers tend to $-\infty$ and that the mass of each component concentrates at its own peaks. We also construct positive multi-peak solutions near prescribed nondegenerate critical points of the respective potentials and compute the leading displacement of each peak. The target frequencies may grow at different rates, provided that $\Lambda_*\to\infty$ and $\log(2+\Lambda^*)=o(\sqrt{\Lambda_*})$, where $\Lambda_*$ and $\Lambda^*$ denote their minimum and maximum. In the $L^2$-noncritical case, the first frequency correction is minus the average potential at the prescribed points. In the $L^2$-critical case, the leading mass term is independent of the frequencies. For each component, we compute the mass correction of order \(\lambda_i^{-2}\) in terms of the sum of the Laplacians of its potential at the prescribed points. When both sums are nonzero, these corrections determine the leading-order asymptotics of the two frequencies and the side from which each self-interaction coefficient approaches its critical value.
\end{abstract}

\maketitle
\enlargethispage{2pt}

\section{Introduction}
We study positive solutions $(u_1,u_2)\in H^1(\mathbb R^N)^2$ of the coupled nonlinear Schr\"odinger system
\begin{equation}\label{eq:system-mu}
\begin{cases}
-\Delta u_1+V_1(x)u_1=a_1u_1^{p-1}+\beta u_1^{\frac{p}{2}-1}u_2^{\frac{p}{2}}+\mu_1u_1,&x\in\mathbb R^N,\\
-\Delta u_2+V_2(x)u_2=a_2u_2^{p-1}+\beta u_2^{\frac{p}{2}-1}u_1^{\frac{p}{2}}+\mu_2u_2,&x\in\mathbb R^N,
\end{cases}
\end{equation}
subject to the two mass constraints
\begin{equation}\label{eq:double-mass}
u_i>0,\qquad \int_{\mathbb R^N}u_i^2\,dx=m_i>0,\qquad i=1,2.
\end{equation}

The potentials are $V_1,V_2$, the self-interaction strengths are $a_1,a_2>0$,
and $\mu_1,\mu_2$ are the Lagrange multipliers for the prescribed masses
$m_1,m_2$. We fix $\beta<0$ throughout. We call a family \emph{segregated}
if the two components concentrate at disjoint sets of points.

The solutions constructed below remain strictly positive throughout $\mathbb R^N$; segregation refers to the separation of their concentration sets.

We work in the Sobolev-subcritical range and assume $p\ge4$ so that the
mixed term has a $C^1$ extension to signed arguments. With $2^*=\infty$ for
$N=1,2$ and $2^*=\frac{2N}{N-2}$ for $N\ge3$, our assumptions are
\begin{equation}\label{eq:p-range}
N\in\{1,2,3\},\qquad 4\le p<2^*.
\end{equation}

For $p=4$ this is a Gross--Pitaevskii system for a two-component
Bose--Einstein condensate. The evolution conserves each component mass
\cite{BaoCai}. For general $p$ in \eqref{eq:p-range},
\eqref{eq:system-mu} is the Euler--Lagrange system for
\begin{equation}\label{eq:energy}
\begin{aligned}
E(u_1,u_2)
={}&\frac12\sum_{i=1}^2\int_{\mathbb R^N}\bigl(|\nabla u_i|^2+V_i u_i^2\bigr)\,dx\\
&-\frac1p\sum_{i=1}^2a_i\int_{\mathbb R^N}|u_i|^p\,dx
-\frac{2\beta}{p}\int_{\mathbb R^N}|u_1|^{\frac{p}{2}}|u_2|^{\frac{p}{2}}\,dx
\end{aligned}
\end{equation}
on the constraint manifold
\begin{equation}\label{eq:additional-1}
\mathcal M_{m_1,m_2}=\left\{(u_1,u_2)\in H^1(\mathbb R^N)^2:
\int_{\mathbb R^N}u_1^2\,dx=m_1,\quad
\int_{\mathbb R^N}u_2^2\,dx=m_2\right\}.
\end{equation}
Each mass constraint has its own Lagrange multiplier.

For scalar equations without a mass constraint, Ambrosetti, Badiale and Cingolani \cite{AmbrosettiBadialeCingolani} constructed semiclassical states near critical points of the potential. Cao, Noussair and Yan \cite{CaoNoussairYan} constructed solutions with peaks near selected critical points of the nonlinear coefficient. We refer to \cite{CaoPengYan} for the reduction method. With prescribed mass, Guo and Seiringer \cite{GuoSeiringer} proved that ground states of a two-dimensional attractive condensate concentrate at a global minimum of the trapping potential as the interaction strength approaches its critical value. Pellacci, Pistoia, Vaira and Verzini \cite{PellacciPistoiaVairaVerzini} constructed normalized solutions concentrating at a nondegenerate critical point of the potential and treated the $L^2$-critical case under a nonvanishing condition. For potentials with nonisolated degenerate critical points, Luo, Peng, Wei and Yan \cite{LuoPengWeiYan} proved existence and local uniqueness of normalized excited states.

For systems with separately prescribed masses, Bartsch and Jeanjean \cite{BartschJeanjean18} studied the autonomous problem by variational methods. Bartsch and Soave \cite{BartschSoave17,BartschSoaveCorrection} used a Pohozaev constraint to study an $L^2$-supercritical competitive system. In the symmetric cubic case, they obtained multiple positive radial solutions for sufficiently strong competition and studied the limit $\beta\to-\infty$ \cite{BartschSoave19}. For positive coupling, Bartsch, Zhong and Zou \cite{BartschZhongZou21} used bifurcation and continuation to construct solutions with prescribed component masses. See also Noris, Tavares and Verzini \cite{NorisTavaresVerzini} for cubic systems with trapping potentials or on bounded domains, and Mederski and Szulkin \cite{MederskiSzulkin} for radial solutions of autonomous multicomponent systems in dimension three. For a two-dimensional nonautonomous system, Tang, Sun and Sun \cite{TangSunSun26} showed that, for sufficiently small values of the spatial perturbation parameter, the numbers of global maxima of the two nonlinear coefficients yield a lower bound on the number of normalized solutions.

Several constructions under a total-mass constraint describe where the peaks of normalized solutions form. Under this constraint, the two components share a common frequency. Guo and Yang \cite{GuoYang23} obtained two-dimensional excited states whose components concentrate at the same set of points. Huang, Pistoia, Troestler and Wang \cite{HuangPistoiaTroestlerWang26} constructed single-peak vector solutions near a nondegenerate critical point of a linear combination of the potentials. Hu, Li, Long and Wang \cite{HuLiLongWang26} proved existence and local uniqueness of normalized multi-peak solutions in dimensions two and three. The two components can also peak at different points under a total-mass constraint, as Guo and Hua \cite{GuoHua25} showed for a slightly $L^2$-supercritical system in dimension two with positive coupling.

For systems without a mass constraint, Peng and Wang \cite{PengWang13} constructed segregated and synchronized solutions in the presence of radial potentials. For sublinear coupling, Guo and Zhang \cite{GuoZhangCV} constructed nonnegative segregated solutions and established a dead-core property. Separation also occurs in the strong-competition limit, where the limiting components have disjoint positivity sets. Related results include the work of Chen and Zou \cite{ChenZou12,ChenZou15} and Clapp and Pistoia \cite{ClappP18} on critical systems, and that of Tavares and Terracini \cite{TavaresTerracini12} on sign-changing solutions with separately prescribed masses.

We prescribe the two component masses separately and keep the competitive coupling fixed. The concentration sets can be chosen independently among nondegenerate critical points of the respective potentials, provided that all prescribed points are distinct. Compared with the total-mass constructions in \cite{HuangPistoiaTroestlerWang26,HuLiLongWang26}, the present problem has two frequencies, whose ratio need not stay bounded or bounded away from zero.

Our estimates must remain uniform as the two frequencies grow at different rates. In particular, when $p=4$, the diagonal term $|\beta|U_{3-i}^2$ in the linearized equation for component $i$ need not be small near the peaks of the other component. We retain this nonnegative term in the diagonal operator and use peak separation to control the off-diagonal terms. This gives a uniform inverse estimate for the projected operator on the admissible frequency rectangles. We then solve the center equations and both mass constraints by a joint degree argument. In the critical case, an explicit second-order mass expansion determines the leading frequencies.

\subsection{Main results}
Whenever $\mu_i<0$, we use the positive frequency
\begin{equation}\label{eq:positive-frequency}
\lambda_i:=-\mu_i>0,\qquad i=1,2.
\end{equation}
The hypotheses on the peak sets and potentials are as follows.

\noindent\textbf{(S)} Let $k_1,k_2$ be positive integers and fix two finite sets
\begin{equation}\label{eq:additional-3}
B_i=\{b_{i,1},\ldots,b_{i,k_i}\},\qquad i=1,2,
\end{equation}
such that the \(k_1+k_2\) points \(b_{i,j}\) are pairwise distinct. Fix \(d_0>0\) so that
\begin{equation}\label{eq:peak-separation}
|b-b'|\ge10d_0\qquad\text{for all distinct }b,b'\in B_1\cup B_2.
\end{equation}

\noindent\textbf{(V1)} For \(i=1,2\), \(V_i\ge0\), \(V_i\in C^1(\mathbb R^N)\), and
\begin{equation}\label{eq:additional-4}
V_i,\ \nabla V_i,\ x\cdot\nabla V_i\in L^\infty(\mathbb R^N).
\end{equation}

\noindent\textbf{(V2)} For \(i=1,2\), \(V_i\) is of class \(C^5\) in a fixed neighborhood of \(B_i\), and its derivatives of orders two through five are bounded there. At each $b_{i,j}$ we assume
\begin{equation}\label{eq:nondegenerate-centers}
\nabla V_i(b_{i,j})=0,\qquad \det D^2V_i(b_{i,j})\ne0.
\end{equation}

\begin{definition}\label{def:concentrating}
Let $a_{1,n},a_{2,n}>0$, and let $(u_{1,n},u_{2,n},\mu_{1,n},\mu_{2,n})$ be solutions of \eqref{eq:system-mu}--\eqref{eq:double-mass}. We call this sequence a segregated \((k_1,k_2)\)-peak concentrating family at \((B_1,B_2)\) if there exist local maximum points \(x_{i,j,n}\to b_{i,j}\) such that
\begin{equation}\label{eq:peak-blow-up}
u_{i,n}(x_{i,j,n})\to+\infty,\qquad i=1,2,\quad j=1,\ldots,k_i,
\end{equation}
and, for every sufficiently small fixed \(\delta>0\),
\begin{equation}\label{eq:exterior-decay}
\|u_{i,n}\|_{L^\infty(\mathbb R^N\setminus\bigcup_{j=1}^{k_i}B_\delta(b_{i,j}))}\to0,
\qquad i=1,2.
\end{equation}

\end{definition}

The exterior region in \eqref{eq:exterior-decay} has infinite measure. The next theorem shows that the mass on this region also tends to zero.

\begin{theorem}\label{thm:necessary}
Assume \eqref{eq:p-range}, \(\beta<0\), (V1), and (S). Let \((u_{1,n},u_{2,n},\mu_{1,n},\mu_{2,n})\) be a segregated concentrating family in the sense of Definition~\ref{def:concentrating}, with arbitrary positive coefficient sequences \(a_{1,n}\) and \(a_{2,n}\). Then
\begin{equation}\label{eq:multiplier-divergence}
\mu_{i,n}\to-\infty,\qquad\lambda_{i,n}\to+\infty,\qquad i=1,2.
\end{equation}
Moreover, for every sufficiently small fixed \(\delta>0\),
\begin{equation}\label{eq:exterior-energy}
\int_{\mathbb R^N\setminus\bigcup_{j=1}^{k_i}B_\delta(b_{i,j})}
\bigl(|\nabla u_{i,n}|^2+\lambda_{i,n}u_{i,n}^2\bigr)\,dx\to0,
\qquad i=1,2.
\end{equation}
In particular,
\begin{equation}\label{eq:mass-concentration}
\int_{\mathbb R^N\setminus\bigcup_{j=1}^{k_i}B_\delta(b_{i,j})}u_{i,n}^2\,dx\to0,
\qquad i=1,2.
\end{equation}
\end{theorem}

Theorem~\ref{thm:necessary} imposes no relation between the two coefficient
sequences. A separate argument is needed when one multiplier remains bounded
while the other tends to $-\infty$. If the multiplier of the first component
remains bounded, the second component decays exponentially near the peaks of
the first. Its contribution to the local Pohozaev identity is then small enough
to bound the $H^1$ norm of the first component. Moser iteration gives a uniform
$L^\infty$ bound, contradicting the divergence of its peak height.

For the construction, put $r=\frac{p}{2}$ and write \eqref{eq:system-mu} as
\begin{equation}\label{eq:frequency-system}
\begin{cases}
-\Delta u_1+(V_1+\lambda_1)u_1=a_1u_1^{p-1}+\beta u_1^{r-1}u_2^r,\\
-\Delta u_2+(V_2+\lambda_2)u_2=a_2u_2^{p-1}+\beta u_2^{r-1}u_1^r.
\end{cases}
\end{equation}
Let \(W\in H^1(\mathbb R^N)\) be the unique positive radial solution of
\begin{equation}\label{eq:ground-state}
-\Delta W+W=W^{p-1}\quad\text{in }\mathbb R^N,
\qquad W(0)=\max_{\mathbb R^N}W.
\end{equation}
We use the uniqueness of $W$ \cite{Kwong} and its nondegeneracy \cite{Weinstein,ChangGustafsonNakanishiTsai}. Set
\begin{equation}\label{eq:ground-state-moments}
M_W:=\int_{\mathbb R^N}W^2\,dx,
\qquad S_W:=\int_{\mathbb R^N}|y|^2W^2(y)\,dy.
\end{equation}
The scaled peak with frequency \(\lambda>0\) and center \(\xi\) is
\begin{equation}\label{eq:scaled-peak}
W_{\lambda,\xi}(x):=\lambda^{\frac{1}{p-2}}W\bigl(\sqrt\lambda(x-\xi)\bigr).
\end{equation}
Its mass satisfies
\begin{equation}\label{eq:mass-scaling}
\int_{\mathbb R^N}W_{\lambda,\xi}^2\,dx=M_W\lambda^q,
\qquad q:=\frac{2}{p-2}-\frac N2.
\end{equation}
The cases \(q<0\), \(q=0\), and \(q>0\) are called the \(L^2\)-supercritical, critical, and subcritical regimes, respectively. Under \eqref{eq:p-range}, the critical case occurs only for \((N,p)=(1,6)\) or \((2,4)\), while the subcritical case occurs only when \(N=1\) and \(4\le p<6\).

Write $\alpha=\frac{2}{p-2}$. For $q\ne0$, the target frequencies are
\begin{equation}\label{eq:noncritical-target}
\Lambda_i(a_i):=\left(\frac{m_i a_i^\alpha}{k_iM_W}\right)^{\frac{1}{q}},
\qquad i=1,2.
\end{equation}
For a target pair, write
\begin{equation}\label{eq:additional-5}
\Lambda_*:=\min\{\Lambda_1,\Lambda_2\},\qquad
\Lambda^*:=\max\{\Lambda_1,\Lambda_2\}.
\end{equation}
In both existence theorems we assume
\begin{equation}\label{eq:frequency-growth}
\Lambda_*\to+\infty,\qquad
\log(2+\Lambda^*)=o(\sqrt{\Lambda_*}).
\end{equation}
This condition permits $\frac{\Lambda_1}{\Lambda_2}$ to tend to zero or infinity.
Across the fixed distance between peaks, the exponential decay at the
slower scale then dominates every fixed power of the faster frequency.

For example, \eqref{eq:frequency-growth} holds for $\Lambda_1=t$ and $\Lambda_2=t^\vartheta$ as $t\to\infty$, for every fixed $\vartheta>0$.

For later use, denote the average potential at the prescribed points by
\[
\overline V_i:=\frac1{k_i}\sum_{j=1}^{k_i}V_i(b_{i,j}),\qquad i=1,2.
\]

\begin{theorem}\label{thm:noncritical}
Assume \eqref{eq:p-range}, \(q\ne0\), \(\beta<0\), (V1), (V2), and (S). Let
\begin{equation}\label{eq:additional-6}
(a_1,a_2)\to
\begin{cases}
(0,0),&q<0,\\
(+\infty,+\infty),&q>0,
\end{cases}
\end{equation}
and suppose that the target frequencies in \eqref{eq:noncritical-target} satisfy \eqref{eq:frequency-growth}. Then \eqref{eq:system-mu}--\eqref{eq:double-mass} has a positive segregated multi-peak solution, with

\begin{equation}\label{eq:frequency-asymptotics}
\lambda_i=\Lambda_i(a_i)-\overline V_i+O\bigl(\Lambda_i(a_i)^{-1}\bigr),
\qquad i=1,2.
\end{equation}

There are strict local maxima \(x_{i,j}\to b_{i,j}\) such that, with \(\omega_{i,j}=\lambda_i+V_i(x_{i,j})\),
\begin{equation}\label{eq:solution-profile}
u_i(x)=a_i^{-\frac{1}{p-2}}\sum_{j=1}^{k_i}\omega_{i,j}^{\frac{1}{p-2}}
W\bigl(\sqrt{\omega_{i,j}}(x-x_{i,j})\bigr)+v_i(x).
\end{equation}
Set \(A_i:=a_i^{-\frac{1}{p-2}}\lambda_i^{\frac{1}{p-2}}\) and \(\lambda_*:=\min\{\lambda_1,\lambda_2\}\). For the fixed constant \(\sigma_0>0\) chosen in the proof, the correction satisfies
\begin{equation}\label{eq:correction-bound}
\begin{split}
|v_i(x)|\le CA_i\lambda_i^{-2}\biggl[&\sum_{j=1}^{k_i}e^{-c\sqrt{\lambda_i}|x-x_{i,j}|}\\
&+e^{-\sigma_0d_0\sqrt{\lambda_*}}\sum_{j=1}^{k_{3-i}}
 e^{-c\sqrt{\lambda_*}|x-x_{3-i,j}|}\biggr].
\end{split}
\end{equation}
Moreover, \(v_i\in H^2(\mathbb R^N)\), and
\begin{equation}\label{eq:center-expansion}
x_{i,j}=b_{i,j}-\frac{S_W}{2NM_W\lambda_i}
\bigl(D^2V_i(b_{i,j})\bigr)^{-1}\nabla\Delta V_i(b_{i,j})+o(\lambda_i^{-1}).
\end{equation}
\end{theorem}

For $q=0$, define
\begin{equation}\label{eq:critical-parameters}
a_{i,*}:=\left(\frac{k_iM_W}{m_i}\right)^{\frac{1}{\alpha}},
\qquad \delta_i(a_i):=k_iM_W-m_i a_i^\alpha,
\end{equation}
and
\begin{equation}\label{eq:critical-coefficient}
\kappa_i:=\frac{S_W}{2N}\sum_{j=1}^{k_i}\Delta V_i(b_{i,j}),\qquad i=1,2.
\end{equation}

\begin{theorem}\label{thm:critical}
Assume \eqref{eq:p-range}, \(q=0\), \(\beta<0\), (V1), (V2), and (S), and suppose that \(\kappa_i\ne0\) for \(i=1,2\). Let \(a_i\to a_{i,*}\) with
\begin{equation}\label{eq:critical-direction}
\delta_i(a_i)\kappa_i>0,\qquad i=1,2,
\end{equation}
and define
\begin{equation}\label{eq:critical-target}
\Lambda_i(a_i):=\left(\frac{\kappa_i}{\delta_i(a_i)}\right)^{\frac{1}{2}},\qquad i=1,2.
\end{equation}
Suppose that this target pair satisfies \eqref{eq:frequency-growth}. Then \eqref{eq:system-mu}--\eqref{eq:double-mass} has a positive segregated multi-peak solution, and
\begin{equation}\label{eq:critical-frequency}
\lambda_i=\Lambda_i(a_i)(1+o(1)),\qquad i=1,2.
\end{equation}
The solution has the profile \eqref{eq:solution-profile} and satisfies \eqref{eq:correction-bound}. The centers $x_{i,j}$ are strict local maxima of $u_i$ and satisfy \eqref{eq:center-expansion}.
\end{theorem}

Condition \eqref{eq:critical-direction} means that $a_i<a_{i,*}$ when $\kappa_i>0$, and $a_i>a_{i,*}$ when $\kappa_i<0$. In particular, nondegenerate local minima give a positive sum of Laplacians, whereas nondegenerate local maxima give a negative sum. For a general collection of critical points, the sum may vanish even though every Hessian is invertible; the assumption $\kappa_i\ne0$ excludes this cancellation.

We first solve the projected problem with the frequencies and peak centers fixed. The peak-center equations and the two mass equations are then solved simultaneously using Brouwer degree. In the $L^2$-critical case, the leading-order mass is independent of the frequency, so a higher-order expansion is required. Testing the localized equation against the scaling derivative of $W$ yields the correction $-\kappa_i/\lambda_i^2$ in the expansion of $a_i^\alpha\int_{\mathbb R^N}u_i^2\,dx$, which determines the leading order of the frequency.

Section~\ref{sec:necessity} proves Theorem~\ref{thm:necessary}. Sections~\ref{sec:approx}--\ref{sec:centers} treat the
projected problem and the center equations. We impose the two mass
constraints in Sections~\ref{sec:noncritical-mass} and~\ref{chap:L2-critical}.

\section{The multipliers and concentration of mass}\label{sec:necessity}

Throughout this section, the solutions form a segregated concentrating
family in the sense of Definition~\ref{def:concentrating}.
Each component \(u_i\) belongs to \(C_{\mathrm{loc}}^2(\mathbb R^N)\).
Indeed, for each fixed \(n\),
\[
-\Delta u_{i,n}
\le (\mu_{i,n})_+u_{i,n}+a_{i,n}u_{i,n}^{p-1}.
\]
Lemma~\ref{lem:moser} gives boundedness of \(u_{i,n}\), and interior
elliptic estimates for the system yield the stated regularity.
We omit the index \(n\) when deriving identities for a single solution.

\subsection{Global identities and spectral comparison}
For a positive solution of \eqref{eq:system-mu}, set
\begin{equation}\label{eq:global-quantities}
K_i:=\int_{\mathbb R^N}|\nabla u_i|^2\,dx,\qquad
P_i:=\int_{\mathbb R^N}u_i^p\,dx,\qquad
I:=\int_{\mathbb R^N}u_1^ru_2^r\,dx,
\end{equation}
and
\begin{equation}\label{eq:global-potentials}
\mathcal V_i:=\int_{\mathbb R^N}(x\cdot\nabla V_i)u_i^2\,dx,\qquad
B:=\sum_{i=1}^2\int_{\mathbb R^N}(V_i-\mu_i)u_i^2\,dx.
\end{equation}
Multiplying the \(i\)th equation by \(u_i\) and integrating gives
\begin{equation}\label{eq:component-nehari}
K_i+\int_{\mathbb R^N}(V_i-\mu_i)u_i^2\,dx=a_iP_i+\beta I.
\end{equation}
Their sum is
\begin{equation}\label{eq:global-nehari}
K_1+K_2+B=a_1P_1+a_2P_2+2\beta I.
\end{equation}
To obtain the global Pohozaev identity, test the equations with
\(\chi_R x\cdot\nabla u_i\), add them and let \(R\to\infty\). This gives
\begin{equation}\label{eq:global-pohozaev}
\frac{N-2}{2}(K_1+K_2)+\frac N2B+\frac12(\mathcal V_1+\mathcal V_2)
=\frac Np(a_1P_1+a_2P_2)+\frac{2N\beta}{p}I.
\end{equation}
The cutoff errors are estimated in Appendix~\ref{app:nehari-pohozaev}.

Subtracting \(\frac{N}{p}\) times \eqref{eq:global-nehari} from \eqref{eq:global-pohozaev} gives
\begin{equation}\label{eq:global-combination}
(1-\gamma)(K_1+K_2)=\gamma B+\frac12(\mathcal V_1+\mathcal V_2),\qquad
\gamma:=\frac{N(p-2)}{2p}\in(0,1).
\end{equation}
The sum of the two coupling terms is a derivative of \(u_1^ru_2^r\).
To estimate a single component we keep the spatial cutoff, since
the separate coupling terms also contain the factor \(x\).

We shall use the following comparison on a bounded domain. For its annular
form, see Luo, Peng, Wei and Yan~\cite[Proposition~2.1]{LuoPengWeiYan}.
\begin{proposition}\label{prop:spectral}
Let \(\Omega\subset\mathbb R^N\) be open and \(c\in C(\Omega)\). Assume that some bounded \(C^2\) domain \(D\Subset\Omega\) satisfies
\begin{equation}\label{eq:spectral-condition}
\inf_{x\in D}c(x)>\lambda_1(D),
\end{equation}
where \(\lambda_1(D)\) is the first Dirichlet eigenvalue of \(-\Delta\) in \(D\). Then
\begin{equation}\label{eq:spectral-equation}
-\Delta w=c(x)w\qquad\text{in }\Omega
\end{equation}
has no positive \(C^2\) solution.
\end{proposition}
\begin{proof}
Let \(\varphi_1>0\) be the first Dirichlet eigenfunction in \(D\). For a positive solution \(w\), Green's identity and the Hopf lemma give
\[
0<\int_D(c-\lambda_1(D))w\varphi_1\,dx
=\int_{\partial D}w\,\partial_\nu\varphi_1\,dS<0 ,
\]
where \(\nu\) is the outward unit normal to \(\partial D\). This is a contradiction.
\end{proof}

\subsection{Local energy estimates and identities}
Fix \(0<d_1<\frac{d_0}{8}\). At each \(b_{i,j}\) choose \(\chi_{i,j}\in C_c^\infty(B_{4d_1}(b_{i,j}))\) such that
\begin{equation}\label{eq:local-cutoffs}
0\le\chi_{i,j}\le1,\qquad \chi_{i,j}=1\ \text{in }B_{3d_1}(b_{i,j}),\qquad
|\nabla\chi_{i,j}|+|D^2\chi_{i,j}|\le C.
\end{equation}
Set
\[
\chi_i:=\sum_{j=1}^{k_i}\chi_{i,j},\qquad
\Omega_i:=\bigcup_{j=1}^{k_i}B_{4d_1}(b_{i,j}).
\]
Define
\begin{equation}\label{eq:local-quantities}
\begin{aligned}
K_i^{\mathrm{loc}}&:=\int_{\mathbb R^N}\chi_i^2|\nabla u_i|^2\,dx,&
B_i^{\mathrm{loc}}&:=\int_{\mathbb R^N}\chi_i^2(V_i-\mu_i)u_i^2\,dx,\\
P_i^{\mathrm{loc}}&:=\int_{\mathbb R^N}\chi_i^2u_i^p\,dx,&
\mathcal V_i^{\mathrm{loc}}&:=\sum_{j=1}^{k_i}\int_{\mathbb R^N}\chi_{i,j}^2
((x-b_{i,j})\cdot\nabla V_i)u_i^2\,dx,
\end{aligned}
\end{equation}
and
\begin{equation}\label{eq:local-mixed}
\begin{aligned}
I_i^{\mathrm{loc}}&:=\int_{\mathbb R^N}\chi_i^2u_i^ru_{3-i}^r\,dx,\\
J_i^{\mathrm{loc}}&:=\sum_{j=1}^{k_i}\int_{\mathbb R^N}\chi_{i,j}^2u_{3-i}^r
(x-b_{i,j})\cdot\nabla(u_i^r)\,dx,\\
Q_i^{\mathrm{loc}}&:=J_i^{\mathrm{loc}}+\frac N2I_i^{\mathrm{loc}}.
\end{aligned}
\end{equation}
\begin{lemma}\label{lem:local-energy}
Suppose that \(\mu_{i,n}\) and \(a_{i,n}\) are bounded. For every fixed compact set \(K\subset\mathbb R^N\setminus B_i\),
\begin{equation}\label{eq:compact-energy}
\int_K(|\nabla u_{i,n}|^2+u_{i,n}^2+a_{i,n}u_{i,n}^p)\,dx\to0.
\end{equation}
For every sufficiently small fixed \(\delta>0\), there exists \(C_\delta>0\), independent of \(n\), such that
\begin{equation}\label{eq:exterior-bound}
\int_{\mathbb R^N\setminus\bigcup_{j=1}^{k_i}B_\delta(b_{i,j})}
(|\nabla u_{i,n}|^2+u_{i,n}^2)\,dx\le C_\delta.
\end{equation}
\end{lemma}
\begin{proof}
Let \(\zeta\) be a smooth compactly supported cutoff with
\(0\le\zeta\le1\). Testing the \(i\)th equation with
\(\zeta^2u_{i,n}\) and using \(V_i\ge0\) and \(\beta<0\), we obtain
\begin{equation}\label{eq:cutoff-energy}
\begin{aligned}
\int_{\mathbb R^N}|\nabla(\zeta u_{i,n})|^2\,dx
&\le
\int_{\mathbb R^N}
\bigl(|\nabla\zeta|^2+(\mu_{i,n})_+\zeta^2\bigr)u_{i,n}^2\,dx\\
&\quad
+a_{i,n}\int_{\mathbb R^N}\zeta^2u_{i,n}^p\,dx,
\end{aligned}
\end{equation}
where \((\mu_{i,n})_+=\max\{\mu_{i,n},0\}\).

To prove \eqref{eq:compact-energy}, choose a bounded open set \(D\)
and a cutoff \(\zeta\) satisfying
\[
\begin{gathered}
K\Subset D\Subset\mathbb R^N\setminus B_i,\qquad
\zeta\in C_c^\infty(D),\qquad
0\le\zeta\le1,\qquad
\zeta\equiv1\quad\text{in }K.
\end{gathered}
\]
By \eqref{eq:exterior-decay}, \(\|u_{i,n}\|_{L^\infty(D)}\to0\).
Since \(\mu_{i,n}\) and \(a_{i,n}\) are bounded,
\eqref{eq:cutoff-energy} gives
\[
\begin{aligned}
\int_K\bigl(|\nabla u_{i,n}|^2+u_{i,n}^2
+a_{i,n}u_{i,n}^p\bigr)\,dx
\le C\int_D\bigl(u_{i,n}^2+a_{i,n}u_{i,n}^p\bigr)\,dx
\longrightarrow0.
\end{aligned}
\]

To prove \eqref{eq:exterior-bound}, set
\[
E_{i,t}:=\mathbb R^N\setminus\bigcup_{j=1}^{k_i}B_t(b_{i,j}).
\]
For a sufficiently small fixed \(\delta>0\), choose a smooth cutoff
\(\zeta\) satisfying
\[
\begin{gathered}
0\le\zeta\le1,\qquad
\zeta=1\quad\text{on }E_{i,\delta},\qquad
\|\nabla\zeta\|_\infty<\infty,\\
\zeta=0\quad\text{in }\bigcup_{j=1}^{k_i}B_{\delta/2}(b_{i,j}).
\end{gathered}
\]
Apply \eqref{eq:cutoff-energy} after inserting an additional cutoff
at infinity. For each fixed \(n\), let its radius tend to infinity,
using \(u_{i,n}\in H^1(\mathbb R^N)\cap L^p(\mathbb R^N)\).
The resulting inequality and the mass constraint give
\[
\begin{aligned}
\int_{E_{i,\delta}}\bigl(|\nabla u_{i,n}|^2+u_{i,n}^2\bigr)\,dx
&\le \int_{\mathbb R^N}|\nabla(\zeta u_{i,n})|^2\,dx+m_i\\
&\le m_i\left(
1+\|\nabla\zeta\|_\infty^2+(\mu_{i,n})_+
+\|a_{i,n}u_{i,n}^{p-2}\|_{L^\infty(E_{i,\delta/2})}
\right).
\end{aligned}
\]
The last term in parentheses is finite for each \(n\) and tends to zero
by \eqref{eq:exterior-decay} and the boundedness of \(a_{i,n}\).
Thus the right-hand side is uniformly bounded in \(n\), which proves
\eqref{eq:exterior-bound}.
\end{proof}
\begin{lemma}\label{lem:local-identities}
Suppose that \(\mu_{i,n}\) and \(a_{i,n}\) are bounded. Then
\begin{equation}\label{eq:local-nehari}
K_{i,n}^{\mathrm{loc}}+B_{i,n}^{\mathrm{loc}}
=a_{i,n}P_{i,n}^{\mathrm{loc}}+\beta I_{i,n}^{\mathrm{loc}}+o(1),
\end{equation}
and
\begin{equation}\label{eq:local-pohozaev}
(1-\gamma)K_{i,n}^{\mathrm{loc}}=\gamma B_{i,n}^{\mathrm{loc}}
+\frac12\mathcal V_{i,n}^{\mathrm{loc}}+\frac{2\beta}{p}Q_{i,n}^{\mathrm{loc}}+o(1).
\end{equation}
\end{lemma}
\begin{proof}
Testing the \(i\)th equation with \(\chi_i^2u_{i,n}\) gives \eqref{eq:local-nehari}, since all terms involving derivatives of the cutoff are \(o(1)\) by Lemma~\ref{lem:local-energy}.

For each \(j\), test the \(i\)th equation with \(\chi_{i,j}^2(x-b_{i,j})\cdot\nabla u_{i,n}\). Integrating by parts and summing over \(j\), we obtain
\[
\frac{N-2}{2}K_{i,n}^{\mathrm{loc}}+\frac N2B_{i,n}^{\mathrm{loc}}
+\frac12\mathcal V_{i,n}^{\mathrm{loc}}
=\frac Np a_{i,n}P_{i,n}^{\mathrm{loc}}-\frac{2\beta}{p}J_{i,n}^{\mathrm{loc}}+o(1).
\]
The scalar cutoff terms are supported on fixed annuli away from \(B_i\), where
\eqref{eq:compact-energy} applies. Appendix~\ref{app:nehari-pohozaev}
gives their bound. Substitution of \eqref{eq:local-nehari} gives \eqref{eq:local-pohozaev}.
\end{proof}

\subsection{Peak heights and decay between the peak sets}
\begin{lemma}\label{lem:maximum}
Fix \(j\in\{1,2\}\) and suppose that \(\inf_n a_{j,n}>0\). Put \(M_{j,n}=\|u_{j,n}\|_\infty\). Then
\begin{equation}\label{eq:maximum-bound}
a_{j,n}M_{j,n}^{p-2}\le C(1+|\mu_{j,n}|).
\end{equation}
\end{lemma}
\begin{proof}
Suppose that \eqref{eq:maximum-bound} fails. Along a subsequence,
\[
\frac{1+|\mu_{j,n}|}{a_{j,n}M_{j,n}^{p-2}}\to0.
\]
By \eqref{eq:peak-blow-up} and \eqref{eq:exterior-decay},
\(M_{j,n}\to\infty\), and for large \(n\) the global maximum is
attained near \(B_j\). Choose a global maximum point \(y_n\) and
pass to a further subsequence so that \(y_n\to b_{j,l}\in B_j\).
Set
\[
\varepsilon_n=(a_{j,n}M_{j,n}^{p-2})^{-\frac12},
\qquad
z_n(y)=\frac{u_{j,n}(y_n+\varepsilon_n y)}{M_{j,n}}.
\]
Then \(\varepsilon_n\to0\), \(0\le z_n\le1\), \(z_n(0)=1\), and the
rescaled equation is
\begin{equation}\label{eq:rescaled-maximum}
\begin{aligned}
-\Delta z_n
+\varepsilon_n^2\bigl(V_j(y_n+\varepsilon_n y)-\mu_{j,n}\bigr)z_n
&=z_n^{p-1}\\
&\quad+\frac{\beta}{a_{j,n}}z_n^{r-1}
\left(\frac{u_{3-j,n}(y_n+\varepsilon_n y)}{M_{j,n}}\right)^r.
\end{aligned}
\end{equation}
By the contradiction assumption, the boundedness of \(V_j\), peak
separation, and \eqref{eq:exterior-decay}, for every fixed \(R>0\) we have
\[
\varepsilon_n^2\bigl(\|V_j\|_\infty+|\mu_{j,n}|\bigr)\to0,
\qquad
\sup_{|y|\le R}u_{3-j,n}(y_n+\varepsilon_n y)\to0.
\]
Together with \(\inf_n a_{j,n}>0\), \(M_{j,n}\to\infty\), and
\(0\le z_n\le1\), the second limit implies that the mixed term in
\eqref{eq:rescaled-maximum} tends to zero locally uniformly.

By interior elliptic estimates and the compact Sobolev embedding, along
a subsequence \(z_n\to z\) in \(C^1_{\mathrm{loc}}(\mathbb R^N)\).
Passing to the limit in \eqref{eq:rescaled-maximum} and applying the
strong maximum principle, we obtain
\[
-\Delta z=z^{p-1}\quad\text{in }\mathbb R^N,
\qquad 0<z\le1,\qquad z(0)=1.
\]

If \(N=3\), then \(p<6\), and the subcritical Liouville theorem of
Gidas and Spruck~\cite{GidasSpruck} rules out such a solution.
If \(N=2\), let \(\bar z(t)\) denote the circular mean of \(z\) on
\(\partial B_t\). By the divergence theorem,
\[
2\pi t\,\bar z'(t)
=-\int_{B_t}z^{p-1}\,dx
\le-\int_{B_1}z^{p-1}\,dx<0,
\qquad t\ge1.
\]
Integrating from \(1\) to \(t\) gives
\[
\bar z(t)\le
\bar z(1)-\frac{\log t}{2\pi}\int_{B_1}z^{p-1}\,dx,
\qquad t\ge1.
\]
The right-hand side tends to \(-\infty\) as \(t\to\infty\),
contradicting \(\bar z(t)>0\). If \(N=1\), then \(z''=-z^{p-1}<0\). Thus \(z\) is strictly concave
on \(\mathbb R\), contradicting its positivity on the whole line. This proves \eqref{eq:maximum-bound}.
\end{proof}

\begin{lemma}\label{lem:away-coefficient}
Fix \(j\in\{1,2\}\) and suppose that \(\lambda_{j,n}:=-\mu_{j,n}\to+\infty\). Let \(E\subset\mathbb R^N\) be a fixed nonempty set with \(\operatorname{dist}(E,B_j)>0\). Then
\begin{equation}\label{eq:away-coefficient}
\sup_{x\in E}\frac{a_{j,n}u_{j,n}^{p-2}(x)}{\lambda_{j,n}}\to0.
\end{equation}
\end{lemma}
\begin{proof}
Set
\[\delta_n:=\sup_{x\in E}u_{j,n}(x),\qquad M_{j,n}:=\|u_{j,n}\|_{L^\infty(\mathbb R^N)}.\]
For a small fixed \(\delta<\operatorname{dist}(E,B_j)\), the set \(E\)
lies in the interior of the exterior region in \eqref{eq:exterior-decay}.
Continuity therefore gives \(\delta_n\to0\), while \eqref{eq:peak-blow-up}
gives \(M_{j,n}\to\infty\).

When \(a_{j,n}\le1\),
\[\frac{a_{j,n}\delta_n^{p-2}}{\lambda_{j,n}}\le\frac{\delta_n^{p-2}}{\lambda_{j,n}}\to0.\]
On any infinite subsequence where \(a_{j,n}>1\), Lemma~\ref{lem:maximum} gives
\[
\frac{a_{j,n}\delta_n^{p-2}}{\lambda_{j,n}}
=\frac{a_{j,n}M_{j,n}^{p-2}}{\lambda_{j,n}}
\left(\frac{\delta_n}{M_{j,n}}\right)^{p-2}
\le C\left(\frac{\delta_n}{M_{j,n}}\right)^{p-2}\to0.
\]
Thus \eqref{eq:away-coefficient} holds.
\end{proof}

\begin{lemma}\label{lem:asymmetric}
Let \(i\in\{1,2\}\) and \(j=3-i\). Suppose that \(a_{i,n}\) is bounded and
\[\mu_{i,n}\text{ is bounded},\qquad\mu_{j,n}\to-\infty.\]
Then
\begin{equation}\label{eq:mixed-vanishing}
I_{i,n}^{\mathrm{loc}}=o(1),\qquad J_{i,n}^{\mathrm{loc}}=o(1),\qquad Q_{i,n}^{\mathrm{loc}}=o(1).
\end{equation}
Consequently,
\begin{equation}\label{eq:local-gradient-bound}
K_{i,n}^{\mathrm{loc}}=O(1).
\end{equation}
\end{lemma}
\begin{proof}
Put \(\lambda_{j,n}=-\mu_{j,n}\). In \eqref{eq:global-combination} the terms
involving the potentials are bounded by (V1) and the mass constraints.
Since \(\mu_{i,n}\) is bounded, it follows that
\begin{equation}\label{eq:total-gradient-growth}
K_{1,n}+K_{2,n}\le C(1+\lambda_{j,n}).
\end{equation}
The closure of \(\widetilde\Omega_i:=\bigcup_\ell B_{6d_1}(b_{i,\ell})\)
is disjoint from \(B_j\). By Lemma~\ref{lem:away-coefficient},
\(a_{j,n}u_{j,n}^{p-2}\le\frac{\lambda_{j,n}}{2}\) there for large \(n\).
Since \(\beta<0\) and \(V_j\ge0\), the \(j\)th equation gives
\begin{equation}\label{eq:comparison-subsolution}
-\Delta u_{j,n}+\frac12\lambda_{j,n}u_{j,n}\le0\quad\text{in }\widetilde\Omega_i.
\end{equation}
Fix \(b=b_{i,\ell}\), put \(R=6d_1\) and choose \(\sigma>0\)
so that \(4\sigma^2R^2<\frac{1}{4}\). On \(B_R(b)\) consider
\[
\Phi_n(x)=\exp[-\sigma\sqrt{\lambda_{j,n}}(R^2-|x-b|^2)].
\]
For large \(n\), \eqref{eq:exterior-decay} gives \(u_{j,n}\le1=\Phi_n\)
on \(\partial B_R(b)\), and
\[
\left(-\Delta+\tfrac12\lambda_{j,n}\right)\Phi_n
\ge\left(\tfrac14\lambda_{j,n}-2N\sigma\sqrt{\lambda_{j,n}}\right)\Phi_n\ge0.
\]
Comparison yields \(u_{j,n}\le\Phi_n\) on the ball. Applying this bound on
\(B_{4d_1}(b)\) at each peak of component \(i\) gives
\begin{equation}\label{eq:other-peak-decay}
\|u_{j,n}\|_{L^\infty(\Omega_i)}\le Ce^{-c\sqrt{\lambda_{j,n}}}.
\end{equation}
By \eqref{eq:total-gradient-growth} and the mass constraints, \(\|u_{i,n}\|_{H^1}\le C(1+\lambda_{j,n})^{\frac{1}{2}}\). By H\"older's inequality,
\[
 |I_{i,n}^{\rm loc}|+|J_{i,n}^{\rm loc}|
 \le C\|u_{j,n}\|_{L^\infty(\Omega_i)}^r
 \left(\|u_{i,n}\|_r^r+
 \|u_{i,n}\|_{2(r-1)}^{r-1}\|\nabla u_{i,n}\|_2\right).
\]
By \eqref{eq:p-range}, both \(r\) and \(2(r-1)\) lie in the range of
the Sobolev embedding for \(H^1\). Thus the preceding estimate and
\eqref{eq:local-mixed} give \eqref{eq:mixed-vanishing}. Then
\eqref{eq:local-pohozaev} and the boundedness of
\(B_{i,n}^{\mathrm{loc}}\) and \(\mathcal V_{i,n}^{\mathrm{loc}}\)
give \eqref{eq:local-gradient-bound}.
\end{proof}

\subsection{Proof of the necessity theorem}
\begin{lemma}\label{lem:moser}
Let \(w_n\ge0\) be weak subsolutions of
\[-\Delta w_n\le C_0w_n+b_nw_n^{p-1}\quad\text{in }\mathbb R^N,\]
where
\[C_0\ge0,\qquad0\le b_n\le C_1,\qquad\|w_n\|_{H^1(\mathbb R^N)}\le H.\]
Under \eqref{eq:p-range}, there exists a constant \(C\), depending only on \(N,p,C_0,C_1,H\), such that
\begin{equation}\label{eq:moser-bound}
\|w_n\|_{L^\infty(\mathbb R^N)}\le C.
\end{equation}
\end{lemma}
\begin{proof}

For \(N=1\), use \(H^1(\mathbb R)\hookrightarrow L^\infty(\mathbb R)\).
For \(N=2,3\), set
\[
s=\frac{p}{p-2}>\frac N2,\qquad s'=\frac p2,\qquad
c_n=C_0+b_nw_n^{p-2}.
\]
Sobolev embedding gives \(\sup_{n,x}\|c_n\|_{L^s(B_2(x))}\le C\).
Take \(\chi=3\) if \(N=3\), and \(\chi=s'+1\) if \(N=2\).
Choose \(t\in(0,1)\) with \(1/s'=1-t+t/\chi\), and set
\(d=2/(1-t)\). Constants depend only on \(N,p,C_0,C_1,H\).

Fix \(x\in\mathbb R^N\) and \(1\le\rho<R\le2\). Choose
\(\eta\in C_c^\infty(B_R(x))\) with \(0\le\eta\le1\),
\(\eta=1\) on \(B_\rho(x)\), and
\(|\nabla\eta|\le C(R-\rho)^{-1}\).
For \(k\ge1\), test with
\(\eta^2w_n\min(w_n,L)^{2k-2}\) and put
\(F=w_n\min(w_n,L)^{k-1}\). Sobolev, H\"older, interpolation,
and Young inequalities give
\[
\begin{aligned}
\|\eta F\|_{2\chi}^2
&\le Ck^2\left[(R-\rho)^{-2}\|F\|_{L^2(B_R(x))}^2
+\|\eta F\|_2^{2(1-t)}\|\eta F\|_{2\chi}^{2t}\right]\\
&\le\tfrac12\|\eta F\|_{2\chi}^2
+Ck^d(R-\rho)^{-2}\|w_n\|_{L^{2k}(B_R(x))}^{2k}.
\end{aligned}
\]
Absorbing and letting \(L\to\infty\), we obtain
\[
\|w_n\|_{L^{2k\chi}(B_\rho(x))}
\le\left[Ck^d(R-\rho)^{-2}\right]^{1/(2k)}
\|w_n\|_{L^{2k}(B_R(x))}.
\]
Starting from \(k=1\), iterate with \(k=\chi^\ell\),
\(R=1+2^{-\ell}\), and \(\rho=1+2^{-\ell-1}\).
The product of constants converges because
\(\sum_{\ell\ge0}(1+\ell)\chi^{-\ell}<\infty\). Hence
\(\|w_n\|_{L^\infty(B_1(x))}\le C\|w_n\|_{L^2(B_2(x))}\le CH\),
uniformly in \(n,x\), which proves \eqref{eq:moser-bound}.

\end{proof}

\begin{proof}[Proof of Theorem~\ref{thm:necessary}]
Suppose one multiplier \(\mu_{i,n}\) is unbounded above. Pass to a
subsequence on which \(\mu_{i,n}\to+\infty\). Choose a ball
\(D\) with \(\overline D\cap(B_1\cup B_2)=\varnothing\). By \eqref{eq:exterior-decay},
\[\delta_n:=\max_{\ell=1,2}\|u_{\ell,n}\|_{L^\infty(D)}\to0.\]
In \(D\), write the \(i\)th equation as \(-\Delta u_{i,n}=c_{i,n}(x)u_{i,n}\), where
\[c_{i,n}:=\mu_{i,n}-V_i+a_{i,n}u_{i,n}^{p-2}
+\beta u_{i,n}^{r-2}u_{3-i,n}^{r}.\]
Since \(a_{i,n}>0\) and \(r\ge2\),
\[\inf_Dc_{i,n}\ge\mu_{i,n}-\|V_i\|_{L^\infty}-|\beta|\delta_n^{p-2}\to+\infty.\]
For large \(n\) this contradicts Proposition~\ref{prop:spectral} applied on \(D\), with \(\Omega=\mathbb R^N\). Hence both multipliers are bounded above.

Now suppose that \(\mu_{i,n}\) stays bounded along a subsequence.
Then \(a_{i,n}\to0\) on that subsequence. Otherwise \(a_{i,n}\) would
be bounded away from zero on a further subsequence, and
\eqref{eq:maximum-bound} would contradict \eqref{eq:peak-blow-up}.
Write \(j=3-i\). Passing to a further subsequence, we may assume that
\(\mu_{j,n}\) stays bounded or tends to \(-\infty\).

When \(\mu_{j,n}\) is bounded, \eqref{eq:global-combination}, the mass constraints \eqref{eq:double-mass}, and (V1) give
\[K_{1,n}+K_{2,n}=O(1).\]
When \(\mu_{j,n}\to-\infty\), Lemma~\ref{lem:asymmetric} gives \(K_{i,n}^{\mathrm{loc}}=O(1)\). Since \(\chi_i=1\) in \(\bigcup_{\ell=1}^{k_i}B_{3d_1}(b_{i,\ell})\), combining this bound with \eqref{eq:exterior-bound}, applied with any sufficiently small \(\delta\le3d_1\), yields \(K_{i,n}=O(1)\). Thus, in both cases,
\[\|u_{i,n}\|_{H^1(\mathbb R^N)}\le C.\]
Since \(\mu_{i,n}\) and \(a_{i,n}\) are bounded, \(\beta<0\), and \(V_i\ge0\),
\[-\Delta u_{i,n}\le Cu_{i,n}+a_{i,n}u_{i,n}^{p-1}.\]
Lemma~\ref{lem:moser} yields \(\|u_{i,n}\|_{L^\infty(\mathbb R^N)}\le C\), contradicting \eqref{eq:peak-blow-up}. Hence \(\mu_{i,n}\) has no bounded subsequence. Together with its upper bound, this implies
\[\mu_{i,n}\to-\infty,\qquad i=1,2.\]

Fix \(i\) and a sufficiently small \(\delta>0\). Write
\[E_{i,t}:=\mathbb R^N\setminus\bigcup_{\ell=1}^{k_i}B_t(b_{i,\ell}),\qquad t>0.\]
Since \(\lambda_{i,n}=-\mu_{i,n}\to+\infty\), Lemma~\ref{lem:away-coefficient} applied to \(E_{i,\frac{\delta}{2}}\) gives
\[\sup_{E_{i,\frac{\delta}{2}}}\frac{a_{i,n}u_{i,n}^{p-2}}{\lambda_{i,n}}\to0.\]
The equation, \(V_i\ge0\), and \(\beta<0\) therefore imply, for all sufficiently large \(n\),
\begin{equation}\label{eq:exterior-subsolution}
-\Delta u_{i,n}+\frac{\lambda_{i,n}}2u_{i,n}\le0
\quad\text{in }\operatorname{int}E_{i,\frac{\delta}{2}}.
\end{equation}
Choose a smooth cutoff \(\eta\) satisfying
\[
0\le\eta\le1,\quad\eta=0\ \text{in }\bigcup_{\ell=1}^{k_i}B_{\frac{3\delta}{4}}(b_{i,\ell}),\quad
\eta=1\ \text{in }E_{i,\delta},\quad|\nabla\eta|\le C_\delta.
\]
Testing \eqref{eq:exterior-subsolution} with \(\eta^2u_{i,n}\), we obtain
\begin{equation}\label{eq:exterior-test}
\int_{\mathbb R^N}|\nabla(\eta u_{i,n})|^2\,dx
+\frac{\lambda_{i,n}}2\int_{\mathbb R^N}\eta^2u_{i,n}^2\,dx
\le\int_{\mathbb R^N}|\nabla\eta|^2u_{i,n}^2\,dx.
\end{equation}
As in Lemma~\ref{lem:local-energy}, a cutoff at infinity justifies the
test. The right-hand side tends to zero by \eqref{eq:exterior-decay},
since \(\nabla\eta\) is supported on fixed bounded annuli away from \(B_i\).
This proves \eqref{eq:exterior-energy} and, since \(\lambda_{i,n}\to\infty\),
\eqref{eq:mass-concentration}.
\end{proof}

\begin{corollary}\label{cor:gradient-divergence}
Under the hypotheses of Theorem~\ref{thm:necessary},
\[\int_{\mathbb R^N}|\nabla u_{i,n}|^2\,dx\to+\infty,\qquad i=1,2.\]
\end{corollary}
\begin{proof}
If \(\|\nabla u_{i,n}\|_2\) were bounded along a subsequence, the mass
constraint and Sobolev embedding would give \(\|u_{i,n}\|_4\le C\) there.
H\"older's inequality would then give
\[
 \int_{\bigcup_j B_\delta(b_{i,j})} u_{i,n}^2
 \le \Big|\bigcup_j B_\delta(b_{i,j})\Big|^{\frac{1}{2}}\|u_{i,n}\|_4^2
 \le C\delta^{\frac{N}{2}}.
\]
By \eqref{eq:mass-concentration} the integral tends to \(m_i\) for
every fixed small \(\delta\). Letting \(\delta\downarrow0\) contradicts \(m_i>0\).
\end{proof}

\section{Approximate solutions and error estimates}\label{sec:approx}
\subsection{Frequency rectangles and approximate profiles}
Fix $0<\eta_0<\frac{1}{10}$. Let $\Lambda_i$ be given by
\eqref{eq:noncritical-target} if $q\ne0$ and by
\eqref{eq:critical-target} if $q=0$. We take the frequencies in
\begin{equation}\mathcal R_a:=\prod_{i=1}^2[(1-\eta_0)\Lambda_i,(1+\eta_0)\Lambda_i].\label{eq:rectangle}\end{equation}
For $\lambda=(\lambda_1,\lambda_2)\in\mathcal R_a$, write
\begin{equation}\label{eq:sec3-2}\lambda_*:=\min\{\lambda_1,\lambda_2\},\qquad \lambda^*:=\max\{\lambda_1,\lambda_2\}.\end{equation}
The frequency growth condition \eqref{eq:frequency-growth} implies, uniformly on $\mathcal R_a$,
\begin{equation}\lambda_*\to+\infty,\qquad \frac{\log(2+\lambda^*)}{\sqrt{\lambda_*}}\to0.\label{eq:growth}\end{equation}
For every $i,j$, set
\begin{equation}\label{eq:sec3-4}z_{i,j}^0:=-\frac{S_W}{2NM_W}\bigl(D^2V_i(b_{i,j})\bigr)^{-1}\nabla\Delta V_i(b_{i,j}).\end{equation}
Write
\[z_{i,j}:=\lambda_i(\xi_{i,j}-b_{i,j}).\]
Choose $R_0>1$ sufficiently large that $|z_{i,j}^0|<\frac{R_0}{2}$ for all $i,j$, and restrict the peak center parameters to
\begin{equation}|z_{i,j}|\le R_0,\quad\text{equivalently}\quad |\xi_{i,j}-b_{i,j}|\le\frac{R_0}{\lambda_i},\quad i=1,2,\quad j=1,\ldots,k_i.\label{eq:centers}\end{equation}
All estimates below are uniform for $\lambda\in\mathcal R_a$ and the
centers in \eqref{eq:centers}, unless stated otherwise. At each center freeze
the potential by taking
\begin{equation}\label{eq:sec3-6}\omega_{i,j}:=\lambda_i+V_i(\xi_{i,j}),\end{equation}
and set
\begin{equation}\label{eq:sec3-7}U_{i,j}(x):=a_i^{-\frac1{p-2}}\omega_{i,j}^{\frac1{p-2}}W\!\left(\sqrt{\omega_{i,j}}(x-\xi_{i,j})\right),\qquad U_i:=\sum_{j=1}^{k_i}U_{i,j}.\end{equation}
Each single-peak function satisfies
\begin{equation}\label{eq:sec3-8}-\Delta U_{i,j}+\omega_{i,j}U_{i,j}=a_iU_{i,j}^{p-1}.\end{equation}
Let
\begin{equation}\label{eq:sec3-9}A_i:=a_i^{-\frac1{p-2}}\lambda_i^{\frac1{p-2}}.\end{equation}
The mass scaling gives constants $c_0,C_0>0$ such that
\begin{equation}c_0\le a_i^{-\frac2{p-2}}\lambda_i^{\frac2{p-2}-\frac N2}\le C_0,\qquad i=1,2.\label{eq:mass-factor}\end{equation}
In the $L^2$-noncritical case, this follows directly from \eqref{eq:noncritical-target} and \eqref{eq:rectangle}. In the $L^2$-critical case, the exponent of $\lambda_i$ in \eqref{eq:mass-factor} is zero, and $a_i$ remains in a compact subset of $(0,+\infty)$. Since
\[A_i^2=\left(a_i^{-\frac2{p-2}}\lambda_i^{\frac2{p-2}-\frac N2}\right)\lambda_i^{\frac{N}{2}},\]
there exists a constant $C>1$ such that
\begin{equation}C^{-1}\lambda_i^{\frac{N}{4}}\le A_i\le C\lambda_i^{\frac{N}{4}},\qquad i=1,2.\label{eq:amplitude}\end{equation}
\subsection{Exponential overlap estimates}
The ground state $W$ satisfies
\begin{equation}\label{eq:sec3-12}|D^\ell W(y)|\le Ce^{-c|y|},\qquad |\ell|\le2.\end{equation}
Consequently,
\begin{equation}|D^\ell U_{i,j}(x)|\le CA_i\lambda_i^{\frac{|\ell|}{2}}e^{-c\sqrt{\lambda_i}|x-\xi_{i,j}|},\qquad |\ell|\le2.\label{eq:peak-decay}\end{equation}
Set
\begin{equation}\label{eq:sec3-14}\sigma:=\begin{cases}1,&p=4,\\ \min\{1,r-2\},&p>4.\end{cases}\end{equation}
Choose $\nu_1$ such that
\begin{equation}0<\nu_1<\min\left\{\frac c4,\frac14\right\},\label{eq:nu}\end{equation}
where $c$ is the decay constant in \eqref{eq:sec3-12}. Next, choose $\sigma_0>0$ sufficiently small that
\begin{equation}32\sigma_0<\nu_1\sigma,\label{eq:sigma0}\end{equation}
and define the exponential decay factor
\begin{equation}\label{eq:sec3-17}\tau_\lambda:=e^{-2\sigma_0d_0\sqrt{\lambda_*}}.\end{equation}
For every fixed $M>0$ and $m\in[0,8]$, uniformly on $\mathcal R_a$,
\begin{equation}(1+\lambda^*)^M e^{-\nu_1\sigma d_0\sqrt{\lambda_*}}=o(\tau_\lambda^m).\label{eq:exp-poly}\end{equation}
Indeed,
\[\frac{(1+\lambda^*)^M e^{-\nu_1\sigma d_0\sqrt{\lambda_*}}}{\tau_\lambda^m}=\exp\!\left[M\log(1+\lambda^*)-(\nu_1\sigma-2m\sigma_0)d_0\sqrt{\lambda_*}\right],\]
and the exponent tends to $-\infty$ by \eqref{eq:growth} and \eqref{eq:sigma0}.
\begin{lemma}\label{lemma:sec3-1}
Let $s,t\ge\sigma$, $1\le h\le\infty$, and let $P$ be a fixed monomial in $a_1,a_2,\lambda_1,\lambda_2$. Then
\begin{equation}|P|\,\|U_1^sU_2^t\|_{L^h(\mathbb R^N)}=o(\tau_\lambda^4)\label{eq:overlap}\end{equation}
uniformly on $\mathcal R_a$ and on the peak center parameter domain. The estimate also holds with $U_1^s U_2^t$ replaced by
$|D^{\gamma_1}U_1|^s|D^{\gamma_2}U_2|^t$ for $|\gamma_1|,|\gamma_2|\le2$.
\end{lemma}
\begin{proof}
For $i=1,2$, set
\[d_i(x):=\operatorname{dist}\bigl(x,\{\xi_{i,1},\ldots,\xi_{i,k_i}\}\bigr).\]
By (S) and \eqref{eq:centers}, all distinct centers are at least $8d_0$ apart. Thus
\begin{equation}\label{eq:sec3-20}d_1(x)+d_2(x)\ge8d_0.\end{equation}
By \eqref{eq:mass-factor}, \eqref{eq:peak-decay} and $s,t\ge\sigma$, the product is bounded by
\[
 |P|U_1^sU_2^t
 \le C(1+\lambda^*)^M e^{-4c\sigma d_0\sqrt{\lambda_*}}
 e^{-(\frac{c\sigma}{2})\sqrt{\lambda_1}d_1(x)}.
\]
The last exponential has uniformly bounded $L^h$ norm, including $h=\infty$.
Estimate \eqref{eq:exp-poly} gives \eqref{eq:overlap}. Derivatives through
order two contribute only further powers of $\lambda^*$ by \eqref{eq:peak-decay}.
\end{proof}
If two peaks belong to the same component but have distinct indices, (S) and \eqref{eq:centers} give the same lower bound $8d_0$ on their center distance. Repeating the preceding argument with $\lambda_1=\lambda_2=\lambda_i$ yields the corresponding exponential smallness estimates for interactions between distinct peaks of the same component.
\subsection{Weighted Banach spaces and the inverse of the positive diagonal operator}
With $\nu_1$ fixed by \eqref{eq:nu}, for $\ell>0$ and $\xi\in\mathbb R^N$, define
\begin{equation}\label{eq:sec3-21}d_{\ell,\xi}(x):=\left(|x-\xi|^2+\ell^{-1}\right)^{\frac{1}{2}}.\end{equation}
Set
\begin{equation}\begin{aligned}w_i(x)&:=\sum_{j=1}^{k_i}e^{-\nu_1\sqrt{\lambda_i}d_{\lambda_i,\xi_{i,j}}(x)},\\ \widetilde w_i(x)&:=\sum_{j=1}^{k_i}e^{-\nu_1\sqrt{\lambda_*}d_{\lambda_*,\xi_{i,j}}(x)},\\ \theta_i(x)&:=\lambda_i^{-2}w_i(x)+\tau_\lambda\widetilde w_{3-i}(x).\end{aligned}\label{eq:weights}\end{equation}
In the term centered at the other component's peaks we use
$\sqrt{\lambda_*}$, the slower of the two decay rates.
By definition,
\[\nabla d_{\ell,\xi}=\frac{x-\xi}{d_{\ell,\xi}},\qquad |\nabla d_{\ell,\xi}|\le1,\qquad \Delta d_{\ell,\xi}=\frac{N-1}{d_{\ell,\xi}}+\frac{\ell^{-1}}{d_{\ell,\xi}^3}\ge0.\]
Thus, for $\ell=\lambda_i$ or $\ell=\lambda_*$,
\[\begin{aligned}\left(-\Delta+\frac{\lambda_i}{4}\right)e^{-\nu_1\sqrt\ell d_{\ell,\xi}}&=\left(\frac{\lambda_i}{4}-\nu_1^2\ell|\nabla d_{\ell,\xi}|^2+\nu_1\sqrt\ell\Delta d_{\ell,\xi}\right)e^{-\nu_1\sqrt\ell d_{\ell,\xi}}\\&\ge\left(\frac14-\nu_1^2\right)\lambda_i e^{-\nu_1\sqrt\ell d_{\ell,\xi}},\end{aligned}\]
where we used $\ell\le\lambda_i$. Summing over the terms in \eqref{eq:weights}, we obtain
\begin{equation}\left(-\Delta+\frac{\lambda_i}{4}\right)\theta_i\ge c\lambda_i\theta_i\quad\text{in }\mathbb R^N,\qquad i=1,2,\label{eq:supersolution}\end{equation}
where $c>0$ is independent of the frequencies and peak center parameters.
\begin{lemma}\label{lemma:sec3-2}
Let $i\ne j$, $a,b\ge\sigma$, and let $M>0$ be fixed. Then, uniformly for $x\in\mathbb R^N$,
\begin{equation}(1+\lambda^*)^M\left(w_i^aw_j^b+w_i^a\widetilde w_j^b+\widetilde w_i^aw_j^b+\widetilde w_i^a\widetilde w_j^b\right)=o(\tau_\lambda^4).\label{eq:weight-products}\end{equation}
Moreover,
\begin{equation}(1+\lambda^*)^M\frac{w_i^aw_j^b}{\widetilde w_j}=o(\tau_\lambda^4),\qquad a\ge\sigma,\quad b\ge1+\sigma.\label{eq:weight-ratio}\end{equation}
\end{lemma}
\begin{proof}
Since \(\lambda_i\ge\lambda_*\),
\[
\sqrt{\lambda_i}\,d_{\lambda_i,\xi}(x)
=\sqrt{1+\lambda_i|x-\xi|^2}
\ge \sqrt{1+\lambda_*|x-\xi|^2}.
\]
Thus \(w_i\le\widetilde w_i\) and, similarly, \(w_j\le\widetilde w_j\).
The sum in parentheses in \eqref{eq:weight-products} is therefore
bounded by \(4\widetilde w_i^a\widetilde w_j^b\).

By (S) and \eqref{eq:centers}, for sufficiently large frequencies,
\[
|x-\xi_{i,\ell}|+|x-\xi_{j,m}|
\ge|\xi_{i,\ell}-\xi_{j,m}|
\ge8d_0.
\]
Using the power estimate for finite sums and \(a,b\ge\sigma\), we obtain
\[
\begin{aligned}
\widetilde w_i^a\widetilde w_j^b
&\le C\sum_{\ell=1}^{k_i}\sum_{m=1}^{k_j}
\exp\!\left[-\nu_1\sigma\sqrt{\lambda_*}
\bigl(|x-\xi_{i,\ell}|+|x-\xi_{j,m}|\bigr)\right]\\
&\le C e^{-8\nu_1\sigma d_0\sqrt{\lambda_*}}.
\end{aligned}
\]
Combining this bound with \eqref{eq:exp-poly} gives
\eqref{eq:weight-products}.

For \eqref{eq:weight-ratio}, the inequality \(w_j\le\widetilde w_j\)
gives
\[
\frac{w_i^aw_j^b}{\widetilde w_j}
\le w_i^aw_j^{b-1}.
\]
Since \(b-1\ge\sigma\), applying \eqref{eq:weight-products} to the
right-hand side proves \eqref{eq:weight-ratio}.
\end{proof}
For continuous pairs $v=(v_1,v_2)$ and $h=(h_1,h_2)$, define
\begin{equation}\label{eq:sec3-26}\|v\|_*:=\max_{i=1,2}\sup_{x\in\mathbb R^N}\frac{|v_i(x)|}{A_i\theta_i(x)},\end{equation}
and
\begin{equation}\label{eq:sec3-27}\|h\|_{**}:=\max_{i=1,2}\sup_{x\in\mathbb R^N}\frac{|h_i(x)|}{\lambda_iA_i\theta_i(x)}.\end{equation}
Let
\begin{equation}\label{eq:sec3-28}\begin{aligned}X_{\lambda,\xi}&:=\{v\in C(\mathbb R^N)^2:\|v\|_*<+\infty\},\\ Y_{\lambda,\xi}&:=\{h\in C(\mathbb R^N)^2:\|h\|_{**}<+\infty\}.\end{aligned}\end{equation}
Since each $\theta_i$ is strictly positive, these are Banach spaces under their respective norms. Moreover, the exponential decay of the weights gives
\[X_{\lambda,\xi}\subset C_0(\mathbb R^N)^2\subset L^\infty(\mathbb R^N)^2,\]
and $\|h\|_{**}<+\infty$ implies $h_i\in L^2(\mathbb R^N)$.
Put $j=3-i$ and define
\begin{equation}\label{eq:sec3-29}\begin{aligned}&Q_i(x):=V_i(x)+\lambda_i-\beta(r-1)U_i^{r-2}(x)U_j^r(x),\\ &D_i\varphi:=-\Delta\varphi+Q_i(x)\varphi.\end{aligned}\end{equation}
Since $\beta<0$ and $V_i\ge0$, we have $Q_i(x)\ge\lambda_i>0$.
\begin{lemma}\label{lemma:sec3-3}
For every continuous function $h_i$ satisfying
\[\sup_{x\in\mathbb R^N}\frac{|h_i(x)|}{\lambda_iA_i\theta_i(x)}<+\infty,\]
the equation
\begin{equation}\label{eq:sec3-30}D_i\varphi_i=h_i\quad\text{in }\mathbb R^N\end{equation}
has a unique solution
\[\varphi_i\in H^1(\mathbb R^N)\cap C^1(\mathbb R^N),\qquad \varphi_i(x)\to0\quad(|x|\to\infty),\]
and
\begin{equation}\sup_{x\in\mathbb R^N}\frac{|\varphi_i(x)|}{A_i\theta_i(x)}\le C\sup_{x\in\mathbb R^N}\frac{|h_i(x)|}{\lambda_iA_i\theta_i(x)}.\label{eq:inverse-bound}\end{equation}
Here $C$ is independent of the parameters. The solution also belongs to $H^2(\mathbb R^N)$.
\end{lemma}
\begin{proof}
Fix the parameters. By the weighted bound on \(h_i\) and the
exponential decay of \(\theta_i\),
\[
h_i\in L^2(\mathbb R^N)\cap L^\infty(\mathbb R^N).
\]
By \eqref{eq:sec3-29}, \(Q_i\in L^\infty(\mathbb R^N)\) and
\(Q_i\ge\lambda_i>0\). The Lax--Milgram theorem gives a unique weak
solution \(\varphi_i\in H^1(\mathbb R^N)\). Since
\[
-\Delta\varphi_i=h_i-Q_i\varphi_i\in L^2(\mathbb R^N),
\]
we have \(\varphi_i\in H^2(\mathbb R^N)\). For \(N\le3\), Sobolev
embedding and interior \(W^{2,s}\) estimates with \(s>N\) give
\(\varphi_i\in C^1(\mathbb R^N)\).

Set
\[
B:=\sup_{x\in\mathbb R^N}
\frac{|h_i(x)|}{\lambda_iA_i\theta_i(x)},
\qquad
\Phi_i:=c^{-1}BA_i\theta_i,
\]
where \(c>0\) is the constant in \eqref{eq:supersolution}. By
\eqref{eq:supersolution} and \(Q_i\ge\lambda_i\),
\[
D_i\Phi_i\ge B\lambda_iA_i\theta_i\ge |h_i|,
\qquad
D_i(\pm\varphi_i-\Phi_i)\le0.
\]
Testing the latter inequalities with
\(q_i^\pm:=(\pm\varphi_i-\Phi_i)_+\in H^1(\mathbb R^N)\) gives
\[
\int_{\mathbb R^N}
\bigl(|\nabla q_i^\pm|^2+\lambda_i(q_i^\pm)^2\bigr)\,dx\le0.
\]
Thus \(q_i^\pm=0\), and hence
\[
|\varphi_i(x)|\le c^{-1}BA_i\theta_i(x).
\]
This proves \eqref{eq:inverse-bound}. Since \(\theta_i(x)\to0\) as
\(|x|\to\infty\), we also have \(\varphi_i(x)\to0\).
\end{proof}
Denote the inverse in Lemma~\ref{lemma:sec3-3} by $T_i:=D_i^{-1}$, and define
\begin{equation}\label{eq:sec3-32}T(h_1,h_2):=(T_1h_1,T_2h_2).\end{equation}
Then $T:Y_{\lambda,\xi}\to X_{\lambda,\xi}$ is a bounded linear operator, with
\[\|Th\|_*\le C\|h\|_{**}.\]
\subsection{The error of the approximate solution}
For $s,t\in\mathbb R$, define
\begin{equation}\label{eq:sec3-33}f(s):=|s|^{p-2}s,\qquad g(s,t):=|s|^{r-2}s|t|^r.\end{equation}
Both maps are continuously differentiable for $p\ge4$. Put
\begin{equation}\label{eq:sec3-34}S_i(u_1,u_2):=-\Delta u_i+(V_i+\lambda_i)u_i-a_if(u_i)-\beta g(u_i,u_{3-i}).\end{equation}
Since $U_i>0$, the residual $R_i:=S_i(U_1,U_2)$ is
\begin{equation}\label{eq:sec3-35}\begin{aligned}R_i={}&\sum_{j=1}^{k_i}(V_i(x)-V_i(\xi_{i,j}))U_{i,j}\\&-a_i\left(U_i^{p-1}-\sum_{j=1}^{k_i}U_{i,j}^{p-1}\right)-\beta U_i^{r-1}U_{3-i}^r.\end{aligned}\end{equation}
\begin{proposition}\label{proposition:sec3-4}
Uniformly on the frequency rectangle and the peak center parameter domain,
\begin{equation}\|R\|_{**}\le C.\label{eq:residual-norm}\end{equation}
More precisely,
\begin{equation}|R_i(x)|\le CA_i\lambda_i^{-1}w_i(x)+o(1)\lambda_iA_i\tau_\lambda\widetilde w_{3-i}(x).\label{eq:residual-pointwise}\end{equation}
\end{proposition}
\begin{proof}
Fix $r_0>0$ sufficiently small that the derivatives of the potentials used below are uniformly bounded near the peak centers by (V2). In $B_{r_0}(\xi_{i,j})$, the center restriction \eqref{eq:centers} and Taylor's formula give
\[|V_i(x)-V_i(\xi_{i,j})|\le C\left(\lambda_i^{-1}|x-\xi_{i,j}|+|x-\xi_{i,j}|^2\right).\]
Setting $y=\sqrt{\lambda_i}(x-\xi_{i,j})$, we obtain from \eqref{eq:peak-decay}
\[|(V_i(x)-V_i(\xi_{i,j}))U_{i,j}(x)|\le CA_i\lambda_i^{-1}(1+|y|^2)e^{-c|y|}.\]
For $x\notin B_{r_0}(\xi_{i,j})$, (V1) and \eqref{eq:peak-decay} yield
\[|(V_i(x)-V_i(\xi_{i,j}))U_{i,j}(x)|\le CA_i e^{-c\sqrt{\lambda_i}|x-\xi_{i,j}|}.\]
Combining this with the estimate inside the ball and \eqref{eq:nu}, and summing over the peaks, gives
\[\left|\sum_{j=1}^{k_i}(V_i(x)-V_i(\xi_{i,j}))U_{i,j}\right|\le CA_i\lambda_i^{-1}w_i.\]
For the interaction error between distinct peaks of the same component, let $d=p-1\ge3$. Then
\[0\le\left(\sum_jU_{i,j}\right)^d-\sum_jU_{i,j}^d\le C\sum_{j\ne\ell}\left(U_{i,j}^{d-1}U_{i,\ell}+U_{i,j}U_{i,\ell}^{d-1}\right).\]
For $j\ne\ell$, use \eqref{eq:peak-decay} on $U_{i,j}^{d-1}U_{i,\ell}$.
One exponential centered at $\xi_{i,j}$ is bounded by $Cw_i$.
The remaining exponentials give $e^{-c'd_0\sqrt{\lambda_i}}$ since
$d-2\ge1$ and the centers are separated. Thus
\[
a_i\left|U_i^{p-1}-\sum_j U_{i,j}^{p-1}\right|
=o(1)A_i\lambda_i^{-1}w_i.
\]
For the mixed term put $j=3-i$. Since $U_h\le CA_hw_h$ and
$\frac{A_i^{r-2}A_j^r}{\lambda_i}$ is bounded by a fixed power of $1+\lambda^*$,
\eqref{eq:weight-ratio} gives
\[
\frac{|\beta|U_i^{r-1}U_j^r}{\lambda_iA_i\tau_\lambda\widetilde w_j}
\le C\frac{A_i^{r-2}A_j^r}{\lambda_i\tau_\lambda}
        \frac{w_i^{r-1}w_j^r}{\widetilde w_j}=o(1).
\]
This proves \eqref{eq:residual-pointwise}. Division by $\lambda_iA_i\theta_i$
gives \eqref{eq:residual-norm}.
\end{proof}

\section{The linearized problem}\label{sec:linear}

We impose orthogonality to the translation modes at each peak. The normalization
below keeps the projections bounded when the two frequencies have different orders.

\subsection{Translation modes and their Gram matrices}

For $l=1,\ldots,N$, let
\begin{equation}
Z_{i,j,l}:=\partial_{x_l}U_{i,j}.
\label{eq:4.1}
\end{equation}
For each fixed component $i$, define the Gram matrix $M_i$, indexed by $(j,l)$ and $(s,m)$, by
\begin{equation}
(M_i)_{(j,l),(s,m)}:=\int_{\R^N}Z_{i,j,l}Z_{i,s,m}\,dx.
\label{eq:4.2}
\end{equation}

\begin{lemma}\label{lem:4.1}
There exist constants $c,C>0$, independent of the parameters, such that, for every coefficient vector $\zeta=(\zeta_{j,l})$,
\begin{equation}
c\lambda_i|\zeta|^2
\le
\sum_{j,l,s,m}(M_i)_{(j,l),(s,m)}\zeta_{j,l}\zeta_{s,m}
\le C\lambda_i|\zeta|^2.
\label{eq:4.3}
\end{equation}
Consequently, $M_i$ is invertible and $\|M_i^{-1}\|\le C\lambda_i^{-1}$.
\end{lemma}

\begin{proof}
Radial symmetry and the change of variables $y=\sqrt{\omega_{i,j}}(x-\xi_{i,j})$ give
\[
 (M_i)_{(j,l),(j,m)}
 =c_W a_i^{-\alpha}\omega_{i,j}^{q+1}\delta_{lm},
 \qquad c_W=\int_{\R^N}|\partial_1W|^2>0.
\]
Since $\frac{\omega_{i,j}}{\lambda_i}\to1$, \eqref{eq:mass-factor} implies
\begin{equation}
 c\lambda_i I_N\le (M_i)^{(j,j)}\le C\lambda_i I_N.
 \label{eq:4.4}
\end{equation}
The off-diagonal blocks are $o(\lambda_i)$ by the overlap estimate for
distinct peaks of one component. There are finitely many blocks, proving
\eqref{eq:4.3} and the inverse bound.
\end{proof}

Define the mass factor
\begin{equation}
s_i:=a_i^{-\alpha}\lambda_i^q.
\label{eq:4.5}
\end{equation}
By \eqref{eq:mass-factor}, $s_i$ is bounded above and below by positive constants uniformly on the prescribed parameter rectangles.

Let
\begin{equation}
\widehat Z_{i,j,l}:=
\sum_{s=1}^{k_i}\sum_{m=1}^N
(M_i^{-1})_{(j,l),(s,m)}Z_{i,s,m}
\label{eq:4.6}
\end{equation}
and set
\begin{equation}
Y_{i,j,l}:=s_i\lambda_i^{-\frac{3}{2}}\widehat Z_{i,j,l}.
\label{eq:4.7}
\end{equation}
For $v=(v_1,v_2)\in X_{\lambda,\xi}$, define
\begin{equation}
\ell_{i,j,l}(v):=
s_i^{-1}\lambda_i^{\frac{3}{2}}
\int_{\R^N}v_iZ_{i,j,l}.
\label{eq:4.8}
\end{equation}
Then
\begin{equation}
\ell_{i,j,l}(Y_{h,s,m}e_h)=\delta_{ih}\delta_{js}\delta_{lm}.
\label{eq:4.9}
\end{equation}

Let $G^0_{lm}=\int_{\R^N}\partial_lW\,\partial_mW$ and
$\psi_l=\sum_m(G^0)^{-1}_{lm}\partial_m W=c_W^{-1}\partial_l W$.
The diagonal blocks of $\frac{M_i}{s_i\lambda_i}$ tend to $G^0$, and the other
blocks tend to zero. Inverting this finite matrix and rescaling $Z$ gives
\begin{equation}
 \lambda_i^2 A_i^{-1}Y_{i,j,l}
 \left(\xi_{i,j}+\frac{y}{\sqrt{\omega_{i,j}}}\right)
 \longrightarrow \psi_l(y)
 \quad\text{locally in }C^1.
 \label{eq:4.10}
\end{equation}

\begin{lemma}\label{lem:4.2}
Uniformly on the prescribed parameter rectangles,
\begin{equation}
|Y_{i,j,l}(x)|+\lambda_i^{-\frac{1}{2}}|\nabla Y_{i,j,l}(x)|
\le CA_i\lambda_i^{-2}w_i(x),
\label{eq:4.11}
\end{equation}
\begin{equation}
\sup_{x\in\R^N}
\frac{|D_iY_{i,j,l}(x)|}{\lambda_iA_i\theta_i(x)}\le C,
\label{eq:4.12}
\end{equation}
and
\begin{equation}
|\ell_{i,j,l}(v)|\le C\|v\|_*.
\label{eq:4.13}
\end{equation}
\end{lemma}

\begin{proof}
By Lemma~\ref{lem:4.1}, \eqref{eq:peak-decay}, and the normalization
in \eqref{eq:4.6}--\eqref{eq:4.7},
\[
\begin{aligned}
|Y_{i,j,l}|+\lambda_i^{-1/2}|\nabla Y_{i,j,l}|
&\le Cs_i\lambda_i^{-3/2}(s_i\lambda_i)^{-1}
A_i\sqrt{\lambda_i}\,w_i\\
&=CA_i\lambda_i^{-2}w_i.
\end{aligned}
\]
This proves \eqref{eq:4.11}.

Set \(h=3-i\). Differentiating the single-peak equation gives
\[
D_i Z_{i,s,m}
=
\Bigl[
V_i-V_i(\xi_{i,s})+(p-1)a_iU_{i,s}^{p-2}
-\beta(r-1)U_i^{r-2}U_h^r
\Bigr]Z_{i,s,m}.
\]
The first two coefficients in brackets are bounded in absolute value
by \(C\lambda_i\). Using \eqref{eq:4.6}--\eqref{eq:4.7},
Lemma~\ref{lem:4.1}, and \eqref{eq:peak-decay}, and applying
\eqref{eq:weight-ratio} to the mixed term, we obtain
\[
\begin{aligned}
|D_iY_{i,j,l}|
&\le CA_i\lambda_i^{-1}w_i
+CA_i^{r-1}A_h^r\lambda_i^{-2}w_i^{r-1}w_h^r\\
&\le C\lambda_iA_i
\bigl(\lambda_i^{-2}w_i+\tau_\lambda\widetilde w_h\bigr)
=C\lambda_iA_i\theta_i.
\end{aligned}
\]
This proves \eqref{eq:4.12}.

Finally, \(\tau_\lambda\le C\lambda_i^{-2}\) and the boundedness of the
weights give
\[
\|v_i\|_{L^\infty}\le CA_i\lambda_i^{-2}\|v\|_*.
\]
By \eqref{eq:peak-decay},
\(\|Z_{i,j,l}\|_{L^1}\le CA_i\lambda_i^{1/2-N/2}\).
Hence, by \eqref{eq:4.8} and \(s_i=A_i^2\lambda_i^{-N/2}\),
\[
|\ell_{i,j,l}(v)|
\le s_i^{-1}\lambda_i^{3/2}
\|v_i\|_{L^\infty}\|Z_{i,j,l}\|_{L^1}
\le C\|v\|_*.
\]
This proves \eqref{eq:4.13}.
\end{proof}

\subsection{The full linear operator}

Linearizing \eqref{eq:sec3-34} at $(U_1,U_2)$, we obtain
\begin{equation}
\begin{aligned}
(Lv)_i={}&-\Delta v_i+(V_i+\lambda_i)v_i-(p-1)a_iU_i^{p-2}v_i\\
&-\beta(r-1)U_i^{r-2}U_j^rv_i
-\beta rU_i^{r-1}U_j^{r-1}v_j,
\qquad j\ne i.
\end{aligned}
\label{eq:4.14}
\end{equation}
Using the operator $D_i$ defined in \eqref{eq:sec3-29}, we write
\begin{equation}
Lv=Dv-Pv,
\label{eq:4.15}
\end{equation}
where
\begin{equation}
(Pv)_i:=(p-1)a_iU_i^{p-2}v_i
+\beta rU_i^{r-1}U_j^{r-1}v_j.
\label{eq:4.16}
\end{equation}

\begin{lemma}\label{lem:4.3}
There exists $\delta_\lambda\to0$, uniformly on the prescribed parameter rectangles, such that
\begin{equation}
\max_{i=1,2}\sup_{x\in\R^N}
\frac{|\beta|rU_i^{r-1}U_j^{r-1}|v_j|}
{\lambda_iA_i\theta_i}
\le\delta_\lambda\|v\|_*.
\label{eq:4.17}
\end{equation}
\end{lemma}

\begin{proof}
By homogeneity, it suffices to assume \(\|v\|_*=1\).
Using \(|v_j|\le A_j\theta_j\) and \(U_h\le CA_hw_h\) for \(h=i,j\),
we obtain
\begin{equation}\label{eq:4.18}
\frac{|\beta|rU_i^{r-1}U_j^{r-1}|v_j|}
     {\lambda_iA_i\theta_i}
\le
C\frac{A_i^{r-2}A_j^r}{\lambda_i}
\frac{w_i^{r-1}w_j^{r-1}\theta_j}{\theta_i}.
\end{equation}
Choose a fixed \(M>0\) such that
\[
\frac{A_i^{r-2}A_j^r}{\lambda_i}\le C(1+\lambda^*)^M.
\]

Expanding \(\theta_j=\lambda_j^{-2}w_j+\tau_\lambda\widetilde w_i\)
and using \(\theta_i\ge\tau_\lambda\widetilde w_j\) for the first term
and \(\theta_i\ge\lambda_i^{-2}w_i\) for the second, we obtain
\begin{equation}\label{eq:4.19}
\frac{w_i^{r-1}w_j^{r-1}\theta_j}{\theta_i}
\le
\frac{\lambda_j^{-2}}{\tau_\lambda}
\frac{w_i^{r-1}w_j^r}{\widetilde w_j}
+\lambda_i^2\tau_\lambda w_i^{r-2}w_j^{r-1}\widetilde w_i.
\end{equation}

For the first term, \eqref{eq:weight-ratio} gives
\[
(1+\lambda^*)^M
\frac{\lambda_j^{-2}}{\tau_\lambda}
\frac{w_i^{r-1}w_j^r}{\widetilde w_j}
=\lambda_j^{-2}o(\tau_\lambda^3)=o(1).
\]
For the second term, \(r\ge2\), the boundedness of the weights,
and \eqref{eq:exp-poly} give
\[
(1+\lambda^*)^M\lambda_i^2\tau_\lambda
w_i^{r-2}w_j^{r-1}\widetilde w_i
\le C(1+\lambda^*)^{M+2}\tau_\lambda=o(1).
\]
These estimates are uniform in \(x\) and the parameters.
Combining them with \eqref{eq:4.18} proves \eqref{eq:4.17}.
\end{proof}

Since $a_i A_i^{p-2}=\lambda_i$ and $U_i\le CA_iw_i$, we have
$a_iU_i^{p-2}|v_i|\le C\lambda_i A_i\theta_i\|v\|_*$. Together with
\eqref{eq:4.17} this gives
\begin{equation}
\|Pv\|_{**}\le C\|v\|_*.
\label{eq:4.20}
\end{equation}
Consequently,
\begin{equation}
K:=TP:X_{\lambda,\xi}\longrightarrow X_{\lambda,\xi}
\label{eq:4.21}
\end{equation}
is bounded.

\begin{lemma}\label{lem:4.4}
For every fixed admissible parameter tuple $(a,\lambda,\xi)$, the operator $K$ defined in \eqref{eq:4.21} is compact.
\end{lemma}

\begin{proof}
In this proof the parameters are fixed. Choose
\begin{equation}
 \nu_1<\nu_+<\min\{2\nu_1,\frac{1}{2}\},
 \label{eq:4.22}
\end{equation}
and denote the weights obtained by replacing $\nu_1$ with $\nu_+$ by
$\theta_i^+$. For $\nu=\nu_1,\nu_+$ the corresponding weight satisfies
\begin{equation}
 c_{a,\lambda,\xi}e^{-\nu\sqrt{\lambda_*}|x|}
 \le \theta_{i,\nu}(x)\le C_{a,\lambda,\xi}e^{-\nu\sqrt{\lambda_*}|x|}.
 \label{eq:4.23}
\end{equation}
Indeed, every exponential in the sum decays at least at this rate, and
$\tau_\lambda\widetilde w_{3-i}$ contains a term with exactly this rate.
From \eqref{eq:peak-decay}, every coefficient of $P$ is bounded by
$C_{a,\lambda,\xi} e^{-\nu_1\sqrt{\lambda_*}|x|}$. Thus, for $\|v\|_*\le1$,
\begin{equation}
 |(Pv)_i|\le C_{a,\lambda,\xi}\lambda_i A_i\theta_i^+.
 \label{eq:4.24}
\end{equation}
The comparison argument for \eqref{eq:inverse-bound} also applies at rate $\nu_+$ and yields
\begin{equation}
 |(Kv)_i|\le C_{a,\lambda,\xi} A_i\theta_i^+.
 \label{eq:4.25}
\end{equation}
Since $D_i(Kv)_i=(Pv)_i$, interior $W^{2,s}$ estimates for $s>N$ give
local precompactness in $C^0$. Moreover, \eqref{eq:4.23}--\eqref{eq:4.25} give
\begin{equation}
 \sup_{\|v\|_*\le1}\sup_{|x|\ge R}
 \frac{|(Kv)_i(x)|}{A_i\theta_i(x)}\longrightarrow0
 \quad(R\to\infty).
 \label{eq:4.26}
\end{equation}
Local compactness and this tail estimate give compactness in $X_{\lambda,\xi}$.
\end{proof}

\subsection{A priori estimates and solvability of the projected linear problem}

Consider the projected linear problem
\begin{equation}
\begin{cases}
Lv=h+\displaystyle\sum_{i=1}^2\sum_{j=1}^{k_i}\sum_{l=1}^N
c_{i,j,l}D_iY_{i,j,l}e_i,\\[1ex]
\ell_{i,j,l}(v)=0,
\quad i=1,2,\ j=1,\ldots,k_i,\ l=1,\ldots,N.
\end{cases}
\label{eq:4.27}
\end{equation}
Here $h\in Y_{\lambda,\xi}$ is given, and $(v,c)\in X_{\lambda,\xi}\times\R^{N(k_1+k_2)}$ is unknown.

Since $L=D-P$, $T=D^{-1}$, and $K=TP$, the equation in \eqref{eq:4.27} is equivalent to
\begin{equation}
v-Kv=Th+
\sum_{i=1}^2\sum_{j=1}^{k_i}\sum_{l=1}^N
c_{i,j,l}Y_{i,j,l}e_i.
\label{eq:4.28}
\end{equation}
Write $Yc=\sum_{i,j,l}c_{i,j,l}Y_{i,j,l}e_i$ and
$\ell(v)=(\ell_{i,j,l}(v))_{i,j,l}$. The augmented system is
\[
 \mathcal F(v,c):=(v-Kv-Yc,\ell(v))=(Th,0).
\]
At $K=0$ the inverse sends $(w,d)$ to $(w+Yc,c)$, where
$c=d-\ell(w)$ by \eqref{eq:4.9}. This inverse is bounded by
\eqref{eq:4.11}--\eqref{eq:4.13}. Lemma~\ref{lem:4.4} therefore shows
that $\mathcal F$ is Fredholm of index zero.

\begin{proposition}\label{prop:4.5}
For sufficiently large target frequencies, every solution of \eqref{eq:4.27} satisfies
\begin{equation}
\|v\|_*+\max_{i,j,l}|c_{i,j,l}|
\le C\|h\|_{**}.
\label{eq:4.29}
\end{equation}
with $C$ independent of the parameters.
\end{proposition}

\begin{proof}
If the estimate fails, choose a sequence of admissible rectangles with
$\lambda_{*,n}\to\infty$ and rescale the solutions so that
\begin{equation}
\|v_n\|_*+\max_{i,j,l}|c_{i,j,l,n}|=1,
\qquad \|h_n\|_{**}\to0.
\label{eq:4.30}
\end{equation}
After extraction, each $c_{i,j,l,n}$ has a limit $c^0_{i,j,l}$.

Fix the $j$th peak of the $i$th component and set
\[
x_{i,j,n}(y):=\xi_{i,j,n}+\frac{y}{\sqrt{\omega_{i,j,n}}}.
\]
Define
\begin{equation}
\varphi_{i,j,n}(y):=
\lambda_{i,n}^2A_{i,n}^{-1}v_{i,n}(x_{i,j,n}(y)).
\label{eq:4.31}
\end{equation}
For each fixed $R>0$, the weighted norm and \eqref{eq:growth} give
\[
\sup_{|y|\le R}|\varphi_{i,j,n}(y)|\le C_R.
\]

Dividing the $i$th component of \eqref{eq:4.27} by $A_{i,n}\lambda_{i,n}^{-2}\omega_{i,j,n}$ and making the change of variables $x=\xi_{i,j,n}+\frac{y}{\sqrt{\omega_{i,j,n}}}$, we obtain
\begin{equation}
\begin{aligned}
-\Delta\varphi_{i,j,n}
&+\frac{V_i(x_{i,j,n}(y))+\lambda_{i,n}}{\omega_{i,j,n}}\varphi_{i,j,n}
-(p-1)\frac{a_{i,n}U_{i,n}^{p-2}(x_{i,j,n}(y))}{\omega_{i,j,n}}\varphi_{i,j,n}
+\mathcal R_{i,j,n}(y)\\
&=\frac{\lambda_{i,n}^2}{A_{i,n}\omega_{i,j,n}}h_{i,n}(x_{i,j,n}(y))
+\frac{\lambda_{i,n}^2}{A_{i,n}\omega_{i,j,n}}
\sum_{s=1}^{k_i}\sum_{l=1}^Nc_{i,s,l,n}D_{i,n}Y_{i,s,l,n}(x_{i,j,n}(y)),
\end{aligned}
\label{eq:4.32}
\end{equation}
where \(\mathcal R_{i,j,n}\) denotes the two rescaled coupling terms.
Peak separation and \eqref{eq:growth} give
\(\mathcal R_{i,j,n}\to0\) in \(L^\infty_{\mathrm{loc}}(\mathbb R^N)\).
Moreover, the two coefficients on the left converge locally uniformly
to \(1\) and \(W^{p-2}\), respectively. Interior elliptic estimates
therefore yield, after passing to a subsequence,
\[
\varphi_{i,j,n}\longrightarrow\varphi_{i,j}
\quad\text{in }C^1_{\mathrm{loc}}(\mathbb R^N).
\]

By the local rescaling limit \eqref{eq:4.10}, only the basis-function
terms corresponding to the current peak survive.
Setting \(L_0=-\Delta+1-(p-1)W^{p-2}\) and passing to the limit
in \eqref{eq:4.32}, we obtain
\begin{equation}
L_0\varphi_{i,j}
=\sum_{l=1}^Nc_{i,j,l}^0(-\Delta+1)\psi_l
\qquad\text{in }\R^N.
\label{eq:4.33}
\end{equation}
By \eqref{eq:4.30}, the definition of \(\theta_{i,n}\), and
\eqref{eq:growth},
\[
\|\varphi_{i,j,n}\|_{L^\infty(\mathbb R^N)}
\le C\bigl(1+\lambda_{i,n}^2\tau_{\lambda_n}\bigr)\le C.
\]
This uniform bound and the exponential decay of \(\partial_m W\)
allow us to pass to the limit in the rescaled orthogonality
conditions by dominated convergence, yielding
\begin{equation}
\int_{\R^N}\varphi_{i,j}\partial_mW=0,
\qquad m=1,\ldots,N.
\label{eq:4.34}
\end{equation}
Test \eqref{eq:4.33} against
$\partial_mW$, using a cutoff and its exponential decay. The left side is
zero. The matrix on the right is $c_W^{-1}$ times the $H^1$ Gram matrix of
the linearly independent functions $\partial_l W$. Thus $c_{i,j,l}^0=0$
for every $l$. Lemma~\ref{lem:appendix-bounded-kernel} and \eqref{eq:4.34} then imply $\varphi_{i,j}=0$.

The same argument applies at every peak. Hence, for every fixed $R>0$,
\begin{equation}
\max_{i,j}\sup_{x\in B_{\frac{R}{\sqrt{\omega_{i,j,n}}}}(\xi_{i,j,n})}
\frac{|v_{i,n}(x)|}{A_{i,n}\theta_{i,n}(x)}\longrightarrow0.
\label{eq:4.35}
\end{equation}
The limits $c_{i,j,l}^0$ found at all the peaks vanish, so that
\begin{equation}
\max_{i,j,l}|c_{i,j,l,n}|\longrightarrow0.
\label{eq:4.36}
\end{equation}

For each $i$, consider the complement of its own balls
$B_{\frac{R}{\sqrt{\omega_{i,j,n}}}}(\xi_{i,j,n})$. By \eqref{eq:peak-decay} we
can fix $R$ large enough that $(p-1)a_{i,n}U_{i,n}^{p-2}\le\frac{\lambda_{i,n}}{8}$
on this complement for both components. Kato's inequality, the nonnegative mixed diagonal term, and Lemma~\ref{lem:4.3} give
\begin{equation}
\begin{aligned}
\left(-\Delta+\frac{\lambda_{i,n}}2\right)|v_{i,n}|
&\le |h_{i,n}|+\delta_{\lambda_n}\lambda_{i,n}A_{i,n}\theta_{i,n}\\
&\quad+\sum_{j=1}^{k_i}\sum_{l=1}^N
|c_{i,j,l,n}|\,|D_{i,n}Y_{i,j,l,n}|.
\end{aligned}
\label{eq:4.37}
\end{equation}
Define
\begin{equation}
\eta_n:=\max\left\{
\|h_n\|_{**},\ \delta_{\lambda_n},\ \max_{i,j,l}|c_{i,j,l,n}|,
\ \max_{i,j}\sup_{x\in\partial B_{\frac{R}{\sqrt{\omega_{i,j,n}}}}(\xi_{i,j,n})}
\frac{|v_{i,n}(x)|}{A_{i,n}\theta_{i,n}(x)}
\right\}.
\label{eq:4.38}
\end{equation}
By \eqref{eq:4.35}--\eqref{eq:4.36},
\begin{equation}
\eta_n\longrightarrow0.
\label{eq:4.39}
\end{equation}
The supersolution inequality \eqref{eq:supersolution} and \eqref{eq:4.12} imply, for a sufficiently large $C$,
\begin{equation}
\begin{aligned}
\left(-\Delta+\frac{\lambda_{i,n}}2\right)
(C\eta_nA_{i,n}\theta_{i,n})
&\ge c_0 C\eta_n\lambda_{i,n}A_{i,n}\theta_{i,n}\\
&\ge |h_{i,n}|+\delta_{\lambda_n}\lambda_{i,n}A_{i,n}\theta_{i,n}
+\sum_{j,l}|c_{i,j,l,n}|\,|D_{i,n}Y_{i,j,l,n}|.
\end{aligned}
\label{eq:4.40}
\end{equation}
Here $c_0>0$ is independent of the parameters. On every inner boundary,
\begin{equation}
|v_{i,n}(x)|\le\eta_nA_{i,n}\theta_{i,n}(x)
\le C\eta_nA_{i,n}\theta_{i,n}(x).
\label{eq:4.41}
\end{equation}
Both functions tend to zero at infinity. The comparison principle therefore gives $|v_{i,n}|\le C\eta_nA_{i,n}\theta_{i,n}$ on this complement. Together with \eqref{eq:4.35}, this yields $\|v_n\|_*\to0$, contradicting \eqref{eq:4.30}. Thus \eqref{eq:4.29} holds.
\end{proof}

\begin{lemma}\label{lem:4.6}
Let $(v,c)$ solve \eqref{eq:4.27}. There exists a fixed exponent $M_g>0$ such that
\begin{equation}
|v_i(x)|+\lambda_i^{-\frac{1}{2}}|\nabla v_i(x)|
\le C(\|v\|_*+\|h\|_{**}+|c|)A_i\theta_i(x),
\quad x\in\bigcup_{j=1}^{k_i}B_{d_0}(\xi_{i,j}),
\label{eq:4.42}
\end{equation}
and
\begin{equation}
|\nabla v_i(x)|
\le C(\|v\|_*+\|h\|_{**}+|c|)(1+\lambda^*)^{M_g}A_i\theta_i(x),
\quad x\in\R^N,
\label{eq:4.43}
\end{equation}
where $|c|:=\max_{i,j,l}|c_{i,j,l}|$. The constants are uniform over the target frequency rectangles and the center parameter domain.
\end{lemma}

\begin{proof}
We first prove \eqref{eq:4.42}. Fix $x\in\bigcup_jB_{d_0}(\xi_{i,j})$, choose a sufficiently small fixed $R_2>0$, and set
\[
B:=B_{\frac{2R_2}{\sqrt{\lambda_i}}}(x),
\qquad
\widetilde v_i(y):=v_i\left(x+\frac{y}{\sqrt{\lambda_i}}\right),
\quad y\in B_{2R_2}(0).
\]
Since $|\nabla\log\theta_i|\le C\sqrt{\lambda_i}$,
\begin{equation}
C^{-1}\theta_i(x)\le\theta_i(z)\le C\theta_i(x),
\qquad z\in B.
\label{eq:4.44}
\end{equation}
On $B$, $a_iU_i^{p-2}\le C\lambda_i$ and
$U_i^{r-2}U_{3-i}^{r}=o(1)$ by separation and \eqref{eq:growth}.
After rescaling, the equation can be written as $-\Delta\widetilde v_i+\widetilde v_i=\widetilde F_i$ in $B_{2R_2}(0)$, where Lemmas~\ref{lem:4.2}--\ref{lem:4.3}, the weighted norms, and \eqref{eq:4.44} give
\begin{equation}
|\widetilde F_i(y)|
\le C(\|v\|_*+\|h\|_{**}+|c|)A_i\theta_i(x),
\qquad y\in B_{2R_2}(0).
\label{eq:4.45}
\end{equation}
The interior $W^{2,s}$ estimate and the Sobolev embedding, with $s>N$, directly yield
\[
|\widetilde v_i(0)|+|\nabla\widetilde v_i(0)|
\le C(\|v\|_*+\|h\|_{**}+|c|)A_i\theta_i(x).
\]
Since $\widetilde v_i(0)=v_i(x)$ and $\nabla\widetilde v_i(0)=\lambda_i^{-\frac{1}{2}}\nabla v_i(x)$, \eqref{eq:4.42} follows.

To prove \eqref{eq:4.43}, choose $M_0>\frac{1}{2}$ and put $r_\lambda=(1+\lambda^*)^{-M_0}$. Since $\sqrt{\lambda^*}r_\lambda\to0$,
\begin{equation}
C^{-1}\theta_i(x)\le\theta_i(z)\le C\theta_i(x),
\qquad z\in B_{2r_\lambda}(x).
\label{eq:4.46}
\end{equation}
Write the $i$th equation as $-\Delta v_i=F_i$. By \eqref{eq:amplitude} and the preceding linear estimates, there exists a fixed $M_1>0$ such that
\begin{equation}
\|F_i\|_{L^\infty(B_{2r_\lambda}(x))}
\le C(1+\lambda^*)^{M_1}(\|v\|_*+\|h\|_{**}+|c|)A_i\theta_i(x),
\label{eq:4.47}
\end{equation}
while
\begin{equation}
\|v_i\|_{L^\infty(B_{2r_\lambda}(x))}
\le C\|v\|_*A_i\theta_i(x).
\label{eq:4.48}
\end{equation}
The local gradient estimate gives
\[
|\nabla v_i(x)|
\le C\left(r_\lambda^{-1}\|v_i\|_{L^\infty(B_{2r_\lambda}(x))}
+r_\lambda\|F_i\|_{L^\infty(B_{2r_\lambda}(x))}\right),
\]
and hence \eqref{eq:4.43}, with $M_g=\max\{M_0,M_1-M_0\}+1$.
\end{proof}

\begin{theorem}\label{thm:4.7}
For all sufficiently large frequency pairs satisfying \eqref{eq:growth}, and for every $h\in Y_{\lambda,\xi}$, problem \eqref{eq:4.27} has a unique solution
\[
(v,c)\in X_{\lambda,\xi}\times\R^{N(k_1+k_2)}.
\]
Moreover,
\begin{equation}
\|v\|_*+\max_{i,j,l}|c_{i,j,l}|
\le C\|h\|_{**}.
\label{eq:4.49}
\end{equation}
The solution belongs to $H^2(\R^N)^2$, and there exists a fixed exponent $M_g>0$ such that
\begin{equation}
|v_i(x)|+\lambda_i^{-\frac{1}{2}}|\nabla v_i(x)|
\le C\|h\|_{**}A_i\theta_i(x),
\quad x\in\bigcup_{j=1}^{k_i}B_{d_0}(\xi_{i,j}),
\label{eq:4.50}
\end{equation}
and
\begin{equation}
|\nabla v_i(x)|
\le C\|h\|_{**}(1+\lambda^*)^{M_g}A_i\theta_i(x),
\quad x\in\R^N.
\label{eq:4.51}
\end{equation}
\end{theorem}

\begin{proof}
Proposition~\ref{prop:4.5} gives injectivity of $\mathcal F$. As its Fredholm
index is zero, it is onto. Applying the proposition to the solution gives
\eqref{eq:4.49}.

The weighted estimate implies $v_i\in L^2(\R^N)$. For each fixed parameter tuple, all coefficients in the linear equation are bounded, and its right-hand side belongs to $L^2(\R^N)$. Standard elliptic regularity therefore gives $v_i\in H^2(\R^N)$. Finally, Lemma~\ref{lem:4.6} and \eqref{eq:4.49} yield \eqref{eq:4.50} and \eqref{eq:4.51}.
\end{proof}

\section{The nonlinear correction}\label{sec:nonlinear}

We use the signed powers $f$ and $g$ of \eqref{eq:sec3-33} until the
projection coefficients vanish. Positivity will follow by testing with
the negative part of each component.

For $v=(v_1,v_2)$, write
\begin{equation}
S(U+v)=R+Lv+N(v). \label{eq:5.1}
\end{equation}
The remainder is
\begin{equation}
\begin{aligned}
N_i(v):={}&-a_i\bigl[f(U_i+v_i)-f(U_i)-f'(U_i)v_i\bigr]\\
&-\beta\bigl[g(U_i+v_i,U_j+v_j)-g(U_i,U_j)
-\partial_1g(U_i,U_j)v_i-\partial_2g(U_i,U_j)v_j\bigr],
\qquad j=3-i.
\end{aligned}
\label{eq:5.2}
\end{equation}

\subsection{Estimates for the remainder}

For $p=4$ put
\begin{equation}
 b_\lambda=\max_{i\ne j}\frac{A_j^2}{\lambda_i\lambda_j^2},
 \label{eq:5.13}
\end{equation}
and put $b_\lambda=0$ for $p>4$. Since $A_j^2\le C\lambda_j^{\frac{N}{2}}$ and
$N\le3$, $b_\lambda\to0$ in either case. With $\sigma$ as in
\eqref{eq:sec3-14}, set
\begin{equation}
 \eta_\lambda=\lambda_*^{-1}+\tau_\lambda^{\frac{\sigma}{2}}+b_\lambda .
 \label{eq:5.14}
\end{equation}

\begin{lemma}\label{lem:5.1}
For each $R>0$ there is a constant $C_R$ such that, if
$\|v\|_*,\|\widetilde v\|_*\le R$, then
\begin{align}
 \|N(v)\|_{**}&\le C_R\eta_\lambda,\label{eq:5.15}\\
 \|N(v)-N(\widetilde v)\|_{**}
   &\le C_R\eta_\lambda\|v-\widetilde v\|_*.\label{eq:5.16}
\end{align}
\end{lemma}

\begin{proof}
It suffices to prove the difference estimate, since $N(0)=0$. Write $d=\|v-\widetilde v\|_*$ and let $z$ lie on the segment between $v$ and
$\widetilde v$. Then
\[
|z_h|\le RA_h\theta_h,\quad |v_h-\widetilde v_h|\le d A_h\theta_h,\quad h=1,2.
\]
The mean value formula applied along this segment reduces the estimate to differences of the first derivatives of $f$ and $g$ at $U+z$ and $U$.

For the scalar term, $|f'(s+h)-f'(s)|\le
C(|s|^{p-3}|h|+|h|^{p-2})$. Thus
\[
 \frac{a_i}{\lambda_i}|f'(U_i+z_i)-f'(U_i)|
 \le C_R(w_i^{p-3}\theta_i+\theta_i^{p-2})
 \le C_R(\lambda_i^{-2}+\tau_\lambda).
\]
Here we used $U_i\le CA_iw_i$, $a_iA_i^{p-2}=\lambda_i$ and boundedness
of the weights.

For the mixed term fix $i$ and put $j=3-i$. Set
\[
 D_s=\partial_1g(U_i+z_i,U_j+z_j)-\partial_1g(U_i,U_j),\qquad
 D_t=\partial_2g(U_i+z_i,U_j+z_j)-\partial_2g(U_i,U_j).
\]
After division by $\lambda_iA_i\theta_i$, the contribution to the
difference estimate is bounded by $|\beta|d$ times
\[
 \frac{|D_s|}{\lambda_i}
       +\frac{A_j\theta_j}{\lambda_iA_i\theta_i}|D_t|.
\]
Suppose first that $r>2$. Write $H(s)=|s|^{r-2}s$. Then
$\partial_1g=(r-1)|s|^{r-2}|t|^r$ and $\partial_2g=rH(s)H(t)$.
For $s\ge0$ and real $h$,
\[
 \big||s+h|^{r-2}-s^{r-2}\big|
 \le C\bigl(s^{r-2-\sigma}|h|^\sigma+|h|^{r-2}\bigr).
\]
For $2<r<3$ this is the H\"older estimate with $\sigma=r-2$.
For $r\ge3$ it follows from the mean value bound, with $\sigma=1$.
Also $|H(s+h)-H(s)|\le C(s^{r-2}|h|+|h|^{r-1})$.
Applying these inequalities to one factor at a time gives, with
$\rho_{ij}=\frac{A_i^{r-2}A_j^r}{\lambda_i}$,
\begin{align*}
 \frac{|D_s|}{\lambda_i}
 &\le C_R\rho_{ij}\Big[
 w_j^r(w_i^{r-2-\sigma}\theta_i^\sigma+\theta_i^{r-2})
 +(w_i^{r-2}+\theta_i^{r-2})(w_j^{r-1}\theta_j+\theta_j^r)\Big],\\
 \frac{A_j\theta_j|D_t|}{\lambda_iA_i\theta_i}
 &\le C_R\rho_{ij}\Big[
 (w_i^{r-2}+\theta_i^{r-2})(w_j^{r-1}\theta_j+\theta_j^r)
 +\lambda_i^2w_i^{r-2}(w_j^{r-2}\theta_j^2+\theta_j^r)\Big].
\end{align*}
For the last term we used $\theta_i\ge\lambda_i^{-2}w_i$. For any
$\gamma>0$, \eqref{eq:weights} gives
$\theta_h^\gamma\le C_\gamma(\lambda_h^{-2\gamma}w_h^\gamma+
\tau_\lambda^\gamma\widetilde w_{3-h}^\gamma)$.
When the second summand is used above, its exponent is at least $\sigma$.
All other products have powers of both $w_i$ and $w_j$ at least $\sigma$.
The latter products satisfy \eqref{eq:weight-products}.
Since $\rho_{ij}$ and $\lambda_i^2\rho_{ij}$ grow at most polynomially in
$1+\lambda^*$, \eqref{eq:growth} bounds both lines by
$C_R\tau_\lambda^{\frac{\sigma}{2}}$.

For $r=2$ we have $g(s,t)=st^2$, so
$D_s=2U_j z_j+z_j^2$ and $D_t=2(U_i z_j+U_j z_i+z_i z_j)$.
Consequently,
\[
 \frac{|D_s|}{\lambda_i}+\frac{A_j\theta_j|D_t|}{\lambda_iA_i\theta_i}
 \le C_R\frac{A_j^2}{\lambda_i}
       \left(w_j\theta_j+\theta_j^2+\frac{w_i\theta_j^2}{\theta_i}\right).
\]
The first two terms on the right are bounded by
$C_R b_\lambda+C_R\tau_\lambda^{\frac{1}{2}}$ upon expanding $\theta_j$.
For the last term use both lower bounds on $\theta_i$ to get
\[
 \frac{w_i\theta_j^2}{\theta_i}\le C\left(
 \frac{\lambda_j^{-4}}{\tau_\lambda}\frac{w_iw_j^2}{\widetilde w_j}
 +\lambda_i^2\tau_\lambda^2\widetilde w_i^2\right).
\]
By \eqref{eq:weight-ratio} and \eqref{eq:growth} this contributes
$O(\tau_\lambda)$ even after multiplication by $\frac{A_j^2}{\lambda_i}$.
Integrating the derivative bounds along the segment proves
\eqref{eq:5.16}.
\end{proof}

\subsection{Solution of the projected equation}

Let $Q_{\lambda,\xi}$ be the solution operator of Theorem~\ref{thm:4.7}.
We seek $(v,c)$ satisfying
\begin{equation}
\begin{cases}
S(U+v)=\displaystyle\sum_{i=1}^2\sum_{j=1}^{k_i}\sum_{l=1}^Nc_{i,j,l}D_iY_{i,j,l}e_i,\\
\ell_{i,j,l}(v)=0,\qquad i=1,2,\quad j=1,\ldots,k_i,\quad l=1,\ldots,N.
\end{cases}
\label{eq:5.25}
\end{equation}
By \eqref{eq:5.1}, this problem is equivalent to $(v,c)=Q_{\lambda,\xi}(-R-N(v))$. Let
\[
\pi_X:X_{\lambda,\xi}\times\R^{N(k_1+k_2)}\longrightarrow X_{\lambda,\xi}
\]
denote the projection onto the first component, and define
\begin{equation}
\mathcal T(v):=\pi_XQ_{\lambda,\xi}(-R-N(v)). \label{eq:5.26}
\end{equation}

\begin{theorem}\label{thm:5.2}
There exists a constant $R_1>0$ such that, for all sufficiently large target frequencies satisfying \eqref{eq:frequency-growth}, problem \eqref{eq:5.25} has a unique solution satisfying
\begin{equation}
v\in X_{\lambda,\xi},\qquad \|v\|_*\le R_1. \label{eq:5.27}
\end{equation}
The corresponding coefficients satisfy $\max_{i,j,l}|c_{i,j,l}|\le C$. Moreover,
\begin{equation}
|v_i(x)|\le CA_i\bigl(\lambda_i^{-2}w_i(x)+\tau_\lambda\widetilde w_{3-i}(x)\bigr), \label{eq:5.28}
\end{equation}
and there exists a fixed exponent $M_g>0$ such that
\begin{equation}
|v_i(x)|+\lambda_i^{-\frac{1}{2}}|\nabla v_i(x)|\le CA_i\theta_i(x),
\qquad x\in\bigcup_{j=1}^{k_i}B_{d_0}(\xi_{i,j}), \label{eq:5.29}
\end{equation}
and
\begin{equation}
|\nabla v_i(x)|\le C(1+\lambda^*)^{M_g}A_i\theta_i(x),\qquad x\in\R^N. \label{eq:5.30}
\end{equation}
The correction also satisfies $v\in H^2(\R^N)^2$ and depends continuously on $(a,\lambda,\xi)$ in every slightly weaker weighted norm.
\end{theorem}

\begin{proof}
By Theorem~\ref{thm:4.7},
\begin{equation}
 \|\pi_XQ_{\lambda,\xi}h\|_*\le C\|h\|_{**}.\label{eq:5.31}
\end{equation}
In particular, $\|\mathcal T(0)\|_*\le C_1$ by \eqref{eq:residual-norm}.
Fix $R_1=2C_1+1$. On the closed ball $\|v\|_*\le R_1$,
\eqref{eq:5.31} and Lemma~\ref{lem:5.1} give
\[
 \|\mathcal T(v)\|_*\le C_1+C_{R_1}\eta_\lambda,\qquad
 \|\mathcal T(v)-\mathcal T(w)\|_*
 \le C_{R_1}\eta_\lambda\|v-w\|_* .
\]
For sufficiently large frequencies, $\mathcal T$ maps the ball into
itself and is a contraction there. At its unique fixed point take
\[
(v,c)=Q_{\lambda,\xi}(-R-N(v)).
\]
This gives \eqref{eq:5.25}.
Theorem~\ref{thm:4.7} bounds $|c|$ and yields the stated regularity and
gradient estimates, since $\|R+N(v)\|_{**}\le C$. Continuity is proved below.
\end{proof}

\begin{lemma}\label{lem:5.3}
Fix $0<\nu_0<\nu_1$. Let $\mathcal K$ be a compact subset of the admissible parameter set on which Theorem~\ref{thm:5.2} applies. Then, for all admissible indices $i,j,l$, the following maps are continuous:
\[
(a,\lambda,\xi)\longmapsto v,\qquad
(a,\lambda,\xi)\longmapsto c_{i,j,l},\qquad
(a,\lambda,\xi)\longmapsto\int_{\R^N}(U_i+v_i)^2.
\]
For $v$, we use the norm with rate $\nu_0$ and weights centered and scaled at the limiting parameters.
\end{lemma}

\begin{proof}
Let
\[
\vartheta_n=(a^{(n)},\lambda^{(n)},\xi^{(n)})\longrightarrow
\vartheta_0=(a^{(0)},\lambda^{(0)},\xi^{(0)})
\]
be a sequence in $\mathcal K$, and denote the corresponding solutions of the projected problem by $(v_n,c_n)$. By Theorem~\ref{thm:5.2},
\begin{equation}
\|v_n\|_{*,\nu_1;\vartheta_n}+\max_{i,j,l}|c_{i,j,l,n}|\le C_{\mathcal K}. \label{eq:5.32}
\end{equation}
By the compactness of the parameter set and \eqref{eq:5.32}, the sequence $(v_n)$ satisfies uniform local $L^\infty$ bounds and uniform exponential decay estimates. The local elliptic estimates then yield, after passing to a subsequence, functions $v$ and coefficients $c$ such that
\[
v_n\longrightarrow v\quad\text{in }C^1_{\mathrm{loc}}(\R^N)^2,
\qquad c_{i,j,l,n}\longrightarrow c_{i,j,l}.
\]
The approximate solutions, the linearized operators, and the projection functions converge locally uniformly with the parameters. Passing to the limit in the projected equations and the orthogonality conditions, we conclude that $(v,c)$ solves problem \eqref{eq:5.25} corresponding to the parameter $\vartheta_0$. Passing to the limit pointwise in $|v_{i,n}|\le R_1A_i^{(n)}\theta_i^{(n)}$ puts $v$ in the same uniqueness ball at $\vartheta_0$. By Theorem~\ref{thm:5.2},
\[
v=v_0,\qquad c_{i,j,l}=c_{i,j,l,0}.
\]
Consequently, the above convergence holds for the entire sequence.

For some $\epsilon>0$ and all $n\ge n_0$,
$\nu_1\sqrt{\lambda_*^{(n)}}>\nu_0\sqrt{\lambda_*^{(0)}}+\epsilon$.
The constants in \eqref{eq:4.23} can be taken uniform near the limiting
parameters, so
\begin{equation}
\lim_{R\to\infty}\sup_{n\ge n_0}\max_{i=1,2}\sup_{|x|\ge R}
\frac{A_i^{(n)}\theta_{i,\nu_1}^{(n)}(x)}{A_i^{(0)}\theta_{i,\nu_0}^{(0)}(x)}=0. \label{eq:5.33}
\end{equation}
Combining \eqref{eq:5.32}, \eqref{eq:5.33}, and the local uniform convergence gives
\begin{equation}
\|v_n-v_0\|_{*,\nu_0;\vartheta_0}\longrightarrow0. \label{eq:5.34}
\end{equation}

Finally, $U_i^{(n)}\to U_i^{(0)}$ and $v_{i,n}\to v_{i,0}$ locally uniformly, while \eqref{eq:5.32} provides an integrable exponential majorant independent of $n$. The dominated convergence theorem therefore yields
\begin{equation}
\int_{\R^N}(U_i^{(n)}+v_{i,n})^2\longrightarrow
\int_{\R^N}(U_i^{(0)}+v_{i,0})^2,
\qquad i=1,2. \label{eq:5.35}
\end{equation}
\end{proof}

\subsection{Positivity}\label{sec:positivity}

Suppose that the coefficients $c_{i,j,l}$ vanish. Then $u=U+v$ solves
the system with the signed nonlinearities in \eqref{eq:sec3-33}.

Define $u_i^-:=\max\{-u_i,0\}$. Testing the $i$-th equation with $-u_i^-$ and using $\beta\le0$, we obtain
\begin{equation}
\int_{\R^N}\bigl(|\nabla u_i^-|^2+(V_i+\lambda_i)(u_i^-)^2\bigr)
\le a_i\int_{\R^N}(u_i^-)^p. \label{eq:5.36}
\end{equation}
Since $u_i=U_i+v_i$ and $U_i>0$, $u_i^-\le|v_i|$. By the pointwise estimate \eqref{eq:5.28} and the identity $a_iA_i^{p-2}=\lambda_i$,
\begin{equation}
\sup_{\R^N}\frac{a_i|v_i|^{p-2}}{\lambda_i}
\le C\bigl(\lambda_i^{-2(p-2)}+\tau_\lambda^{p-2}\bigr)=o(1). \label{eq:5.37}
\end{equation}
Hence the right-hand side of \eqref{eq:5.36} is at most
$(\frac{\lambda_i}{2})\|u_i^-\|_2^2$, giving $u_i^-=0$.
At each center, $U_i(\xi_{i,j})\ge A_iW(0)$ and
$|v_i(\xi_{i,j})|=o(A_i)$. The component $u_i$ is therefore nonzero.
It satisfies
\[
 -\Delta u_i+(V_i+\lambda_i+|\beta|u_i^{r-2}|u_{3-i}|^r)u_i
 =a_i u_i^{p-1}\ge0,
\]
where $u_i^{r-2}$ denotes $1$ when $r=2$. The coefficient is locally
bounded and nonnegative, so the strong maximum principle gives $u_i>0$.

\section{The center equations}\label{sec:centers}

To determine the centers we test the projected equation against the
translations of each peak. For the solution $(v,c)$ of
Theorem~\ref{thm:5.2}, put
\begin{equation}
G_{i,j,l}(a,\lambda,\xi)
:=\int_{\R^N}S_i(U+v)Z_{i,j,l}.
\label{eq:6.1}
\end{equation}
This integral is absolutely convergent. By \eqref{eq:5.25},
\begin{equation}
G_{i,j,l}
=\sum_{s,m}c_{i,s,m}\int_{\R^N}D_iY_{i,s,m}Z_{i,j,l}.
\label{eq:6.2}
\end{equation}
By Lemma~\ref{lem:5.3}, $G$ is continuous with respect to all the finite-dimensional
parameters.

\begin{lemma}\label{lemma:sec6-1}
Define the $N\times N$ matrix $H^0$ by
\[
H^0_{lm}:=\int_{\R^N}
\bigl(\nabla\psi_m\cdot\nabla\partial_lW+\psi_m\partial_lW\bigr).
\]
Then $H^0$ is positive definite and, uniformly on the parameter rectangles,
\begin{equation}
\sqrt{\lambda_i}\int_{\R^N}D_iY_{i,s,m}Z_{i,j,l}
=s_iH^0_{lm}\delta_{js}+o(1).
\label{eq:6.3}
\end{equation}
Let $\mathcal B_i$ be the matrix whose entries are
\[
(\mathcal B_i)_{(j,l),(s,m)}
:=\int_{\R^N}D_iY_{i,s,m}Z_{i,j,l}.
\]
Then, for all sufficiently large frequencies, $\mathcal B_i$ is invertible and
\begin{equation}
\lVert \mathcal B_i^{-1}\rVert\leq C\frac{\sqrt{\lambda_i}}{s_i}.
\label{eq:6.4}
\end{equation}
Consequently,
\begin{equation}
G=0\quad\Longleftrightarrow\quad c=0.
\label{eq:6.5}
\end{equation}
\end{lemma}

\begin{proof}
Since $\psi_m=c_W^{-1}\partial_m W$, the matrix $H^0$ is $c_W^{-1}$ times the $H^1$ Gram matrix of the linearly independent functions $\partial_mW$, and is therefore positive definite. Using the local scaling limits in \eqref{eq:4.10}, the identity
\[
A_i^2\lambda_i^{-N/2}=s_i,
\]
and the exponential smallness of the interactions between distinct
peaks, we obtain
\[
\sqrt{\lambda_i}\int_{\mathbb R^N}
D_iY_{i,s,m}Z_{i,j,l}\,dx
=s_iH^0_{lm}\delta_{js}+o(1),
\]
which proves \eqref{eq:6.3}. Consequently,
\[
\mathcal B_i
=\frac{s_i}{\sqrt{\lambda_i}}
\operatorname{diag}(H^0,\ldots,H^0)
+o(\lambda_i^{-1/2}).
\]
Since \(H^0\) is positive definite and \(s_i\) is uniformly bounded
away from zero, \(\mathcal B_i\) is invertible for all sufficiently
large frequencies and
\[
\|\mathcal B_i^{-1}\|\le C\frac{\sqrt{\lambda_i}}{s_i}.
\]
This proves \eqref{eq:6.4}, while \eqref{eq:6.5} follows from
\eqref{eq:6.2}.
\end{proof}

\subsection{Expansion of the center equations}

\begin{proposition}\label{proposition:sec6-2}
For center configurations satisfying
\[
|\xi_{i,j}-b_{i,j}|\leq\frac{R_0}{\lambda_i},
\]
we have
\begin{equation}
\begin{aligned}
G_{i,j,l}
={}&-\frac12a_i^{-\alpha}\omega_{i,j}^{q}
\left[
M_W\partial_lV_i(\xi_{i,j})
+\frac{S_W}{2N\omega_{i,j}}\partial_l\Delta V_i(\xi_{i,j})
\right] \\
&\quad+o\bigl(a_i^{-\alpha}\lambda_i^{q-1}\bigr).
\end{aligned}
\label{eq:6.6}
\end{equation}
The remainder is uniform with respect to the other frequency and center parameters.
\end{proposition}

\begin{proof}
Since \(\nabla V_i(b_{i,j})=0\),
\[
\partial_lV_i(\xi_{i,j})=O(\lambda_i^{-1}).
\]
By \eqref{eq:5.1},
\[
G_{i,j,l}
=\int_{\mathbb R^N}R_iZ_{i,j,l}\,dx
+\int_{\mathbb R^N}(Lv)_iZ_{i,j,l}\,dx
+\int_{\mathbb R^N}N_i(v)Z_{i,j,l}\,dx.
\]
The exponential separation of the peaks gives
\[
\int_{\mathbb R^N}R_iZ_{i,j,l}\,dx
=\int_{\mathbb R^N}
\bigl(V_i(x)-V_i(\xi_{i,j})\bigr)U_{i,j}Z_{i,j,l}\,dx
+o(s_i\lambda_i^{-1}).
\]
Since \(Z_{i,j,l}=\partial_lU_{i,j}\),
\begin{equation}\label{eq:6.7}
\int_{\mathbb R^N}
\bigl(V_i(x)-V_i(\xi_{i,j})\bigr)U_{i,j}Z_{i,j,l}\,dx
=-\frac12\int_{\mathbb R^N}\partial_lV_i(x)U_{i,j}^2(x)\,dx.
\end{equation}
After the change of variables
\[
y=\sqrt{\omega_{i,j}}(x-\xi_{i,j}),
\]
Taylor expansion and radial symmetry yield
\begin{equation}\label{eq:6.8}
\begin{aligned}
\int_{\mathbb R^N}\partial_lV_i(x)U_{i,j}^2(x)\,dx
&=a_i^{-\alpha}\omega_{i,j}^q
\left[
M_W\partial_lV_i(\xi_{i,j})
+\frac{S_W}{2N\omega_{i,j}}\partial_l\Delta V_i(\xi_{i,j})
\right]\\
&\quad+o\bigl(a_i^{-\alpha}\lambda_i^{q-1}\bigr).
\end{aligned}
\end{equation}
Therefore,
\[
\begin{aligned}
\int_{\mathbb R^N}R_iZ_{i,j,l}\,dx
&=-\frac12a_i^{-\alpha}\omega_{i,j}^q
\left[
M_W\partial_lV_i(\xi_{i,j})
+\frac{S_W}{2N\omega_{i,j}}\partial_l\Delta V_i(\xi_{i,j})
\right]\\
&\quad+o\bigl(a_i^{-\alpha}\lambda_i^{q-1}\bigr).
\end{aligned}
\]

We next justify the estimates for the two correction terms. Fix the
indices $i,j,l$, put $h=3-i$, and abbreviate
$U^0=U_{i,j}$, $Z=Z_{i,j,l}$, and $\xi=\xi_{i,j}$.
Choose a fixed $0<\rho<d_0$ such that (V2) holds on
$B_{2\rho}(\xi)$ for all admissible centers, and set
$B=B_\rho(\xi)$. In this argument $E_\lambda$ denotes a nonnegative
quantity bounded by
\[
 C(1+\lambda^*)^M e^{-\gamma\sqrt{\lambda_*}},
\]
where $C,M,\gamma>0$ are independent of the parameters and may change
from one occurrence to the next. By \eqref{eq:growth} and the upper
and lower bounds on $s_i$, for every fixed $K>0$ we have
\[
 E_\lambda=o(s_i\lambda_i^{-K})
 \qquad (i=1,2).
\]
All amplitudes and coefficients appearing below are bounded by a fixed
power of $1+\lambda^*$, by \eqref{eq:mass-factor} and
\eqref{eq:amplitude}.

Peak separation, \eqref{eq:peak-decay}, and \eqref{eq:5.28} imply,
uniformly on $B$, with $y=\sqrt{\lambda_i}(x-\xi)$,
\[
\begin{aligned}
 |U_i-U^0|+|U_h|+|v_h|&\le E_\lambda,\\
 |U^0(x)|&\le CA_i e^{-c|y|},\\
 |Z(x)|&\le CA_i\sqrt{\lambda_i}e^{-c|y|},\\
 |v_i(x)|&\le CA_i\lambda_i^{-2}e^{-\nu_1|y|}+E_\lambda.
\end{aligned}
\]
Here the term $\tau_\lambda\widetilde w_i$ in the estimate for $v_h$
is included in $E_\lambda$. On the complement of $B$,
\[
 \int_{\mathbb R^N\setminus B}|Z|\,dx\le E_\lambda.
\]
Consequently, the integral against $Z$ of any of the algebraic
expressions below over this complement is also bounded by
$E_\lambda$.

Since $v_i\in H^2(\mathbb R^N)$, integration by parts and the equation
for $U^0$ give the exact identity
\[
\begin{aligned}
 \int_{\mathbb R^N}(Lv)_i Z\,dx
 ={}&\int_{\mathbb R^N}v_i
 \Bigl[V_i-V_i(\xi)
 +(p-1)a_i\bigl((U^0)^{p-2}-U_i^{p-2}\bigr)
 -\beta(r-1)U_i^{r-2}U_h^r\Bigr]Z\,dx\\
 &-\beta r\int_{\mathbb R^N}
 U_i^{r-1}U_h^{r-1}v_h Z\,dx.
\end{aligned}
\]
All terms on the right other than the potential difference are bounded
by $E_\lambda$: on $B$ this follows from the preceding local bounds,
and on its complement from the decay of $Z$. This also covers $r=2$,
when $U_i^{r-2}=1$; no global smallness of the diagonal coupling
coefficient is needed. On $B$, (V2) and
$|\xi-b_{i,j}|\le R_0/\lambda_i$ yield
\[
 |V_i(x)-V_i(\xi)|
 \le C\bigl(\lambda_i^{-1}|x-\xi|+|x-\xi|^2\bigr)
 =C\bigl(\lambda_i^{-3/2}|y|+\lambda_i^{-1}|y|^2\bigr).
\]
It follows, using $A_i^2\lambda_i^{-N/2}=s_i$, that
\[
\begin{aligned}
 \int_B|v_i|\,|V_i-V_i(\xi)|\,|Z|\,dx
 &\le Cs_i\lambda_i^{-3/2}
 \int_{\mathbb R^N}
 \bigl(\lambda_i^{-3/2}|y|+\lambda_i^{-1}|y|^2\bigr)
 e^{-c_1|y|}\,dy+E_\lambda\\
 &\le Cs_i\lambda_i^{-5/2}+E_\lambda,
\end{aligned}
\]
where $c_1>0$ is fixed. Thus
\begin{equation}\label{eq:6.9}
 \left|\int_{\mathbb R^N}(Lv)_iZ_{i,j,l}\,dx\right|
 \le Cs_i\lambda_i^{-5/2}+o(s_i\lambda_i^{-3})
 =o(s_i\lambda_i^{-1}).
\end{equation}

For the scalar part of $N_i(v)$ we use, since $p\ge4$,
\[
 |f(s+t)-f(s)-f'(s)t|
 \le C\bigl(s^{p-3}|t|^2+|t|^{p-1}\bigr),
 \qquad s\ge0,\quad t\in\mathbb R.
\]
The same local bounds and $a_iA_i^{p-2}=\lambda_i$ therefore give
\[
\begin{aligned}
 a_i\int_{\mathbb R^N}
 \bigl(U_i^{p-3}|v_i|^2+|v_i|^{p-1}\bigr)|Z|\,dx
 &\le Cs_i\bigl(\lambda_i^{-5/2}
       +\lambda_i^{3/2-2(p-1)}\bigr)+E_\lambda\\
 &\le Cs_i\lambda_i^{-5/2}+E_\lambda.
\end{aligned}
\]
For the mixed part, recall that
\[
 g(s,t)=|s|^{r-2}s|t|^r,\qquad
 \partial_1g(s,t)=(r-1)|s|^{r-2}|t|^r,\qquad
 \partial_2g(s,t)=r|s|^{r-2}s|t|^{r-2}t.
\]
On $B$, each of
\[
 g(U_i+v_i,U_h+v_h),\quad g(U_i,U_h),\quad
 \partial_1g(U_i,U_h)v_i,\quad
 \partial_2g(U_i,U_h)v_h
\]
is bounded by $E_\lambda$, because $U_h$ and $v_h$ are exponentially
small there and $r\ge2$. Their integrals against $Z$ over the
complement of $B$ are bounded in the same way by the decay of $Z$.
Hence the mixed remainder contributes only $E_\lambda$, including
when $p=4$. Combining the scalar and mixed estimates proves
\begin{equation}\label{eq:6.10}
 \left|\int_{\mathbb R^N}N_i(v)Z_{i,j,l}\,dx\right|
 \le Cs_i\lambda_i^{-5/2}+o(s_i\lambda_i^{-3})
 =o(s_i\lambda_i^{-1}).
\end{equation}
Every estimate is uniform in both frequencies and all admissible
centers. Since $s_i\lambda_i^{-1}=a_i^{-\alpha}\lambda_i^{q-1}$,
combining \eqref{eq:6.7}--\eqref{eq:6.10} proves
\eqref{eq:6.6}.

\end{proof}

Set
\begin{equation}
\xi_{i,j}=b_{i,j}+\lambda_i^{-1}z_{i,j},
\qquad |z_{i,j}|\leq R_0,
\label{eq:6.11}
\end{equation}
and write
\[
H_{i,j}:=D^2V_i(b_{i,j}),
\qquad
\gamma_W:=\frac{S_W}{2NM_W}.
\]
By Proposition~\ref{proposition:sec6-2}, the limiting center equation is
\[
H_{i,j}z_{i,j}+\gamma_W\nabla\Delta V_i(b_{i,j})=0,
\]
whose unique solution is
\[
z^0_{i,j}=-\gamma_WH_{i,j}^{-1}\nabla\Delta V_i(b_{i,j}).
\]
Choose $R_c>0$ such that
\[
\overline{B_{R_c}(z^0_{i,j})}\subset B_{R_0}(0),
\]
and set
\begin{equation}
\cZ:=\prod_{i=1}^{2}\prod_{j=1}^{k_i}B_{R_c}(z^0_{i,j}).
\label{eq:6.12}
\end{equation}
Define
\begin{equation}
\cG_{i,j}(z)
:=-\frac{2\lambda_i a_i^\alpha}{M_W\omega_{i,j}^{q}}
\bigl(G_{i,j,1},\ldots,G_{i,j,N}\bigr).
\label{eq:6.13}
\end{equation}
Then, uniformly on $\overline{\cZ}$,
\begin{equation}
\cG_{i,j}(z)=H_{i,j}(z_{i,j}-z^0_{i,j})+o(1).
\label{eq:6.14}
\end{equation}
Since all $H_{i,j}$ are invertible, \eqref{eq:6.14} and the homotopy invariance of the
degree give, for sufficiently large target frequencies,
\begin{equation}
\begin{aligned}
\deg(\cG,\cZ,0)
&=\operatorname{sgn}\prod_{i=1}^{2}\prod_{j=1}^{k_i}\det H_{i,j} \\
&=\operatorname{sgn}\prod_{i=1}^{2}\prod_{j=1}^{k_i}
\det D^2V_i(b_{i,j})\neq0.
\end{aligned}
\label{eq:6.15}
\end{equation}
For each fixed frequency pair there is therefore a zero of $G$ in these
coordinates. Every such zero satisfies
\[
z_{i,j}=z^0_{i,j}+o(1),
\qquad
\xi_{i,j}=b_{i,j}+\lambda_i^{-1}z^0_{i,j}+o(\lambda_i^{-1}).
\]

\subsection{Decay away from the peaks}

\begin{proposition}\label{proposition:sec6-3}
For the projected solution $u=U+v$, with $M_g$ as in Theorem~\ref{thm:4.7}, we have
\begin{equation}
|u_i(x)|\leq CA_i\left[
\sum_{j=1}^{k_i}e^{-c\sqrt{\lambda_i}|x-\xi_{i,j}|}
+\tau_\lambda\sum_{j=1}^{k_{3-i}}
e^{-c\sqrt{\lambda_*}|x-\xi_{3-i,j}|}
\right],
\label{eq:6.16}
\end{equation}
and
\begin{equation}
|\nabla u_i(x)|\leq CA_i\left[
\sqrt{\lambda_i}\sum_{j=1}^{k_i}
e^{-c\sqrt{\lambda_i}|x-\xi_{i,j}|}
+(1+\lambda^*)^{M_g}\theta_i(x)
\right].
\label{eq:6.17}
\end{equation}
Moreover, for every fixed $M>0$,
\begin{equation}
\int_{\R^N}|u_1|^r|u_2|^r
+\sum_{i=1}^{2}\int_{\R^N}|u_i|^{r-1}|u_{3-i}|^r|\nabla u_i|
\leq C_M(1+\lambda^*)^{-M}.
\label{eq:6.18}
\end{equation}
\end{proposition}

\begin{proof}
Estimates \eqref{eq:6.16} and \eqref{eq:6.17} follow from the corresponding bounds for the
principal profiles and from \eqref{eq:5.28} and \eqref{eq:5.30}. In particular, there exists a
fixed exponent $M_0>0$ such that
\[
|u_i|+(1+\lambda^*)^{-M_0}|\nabla u_i|
\leq CA_i(w_i+\tau_\lambda\widetilde w_{3-i}).
\]
Using $(x+y)^r\leq C_r(x^r+y^r)$, the integrands in \eqref{eq:6.18} are bounded by a
fixed algebraic factor in the frequencies times
\[
(w_1+\tau_\lambda\widetilde w_2)^r
(w_2+\tau_\lambda\widetilde w_1)^r
\leq C\bigl(
w_1^rw_2^r+\tau_\lambda^rw_1^r\widetilde w_1^r
+\tau_\lambda^r\widetilde w_2^rw_2^r
+\tau_\lambda^{2r}\widetilde w_2^r\widetilde w_1^r
\bigr).
\]
For the first and fourth products, apply \eqref{eq:weight-products} to half of each power and integrate the remaining exponentials. Their integrals are uniformly bounded since $\lambda_*\ge1$. The middle two products also have bounded integrals and carry the factor $\tau_\lambda^r$. Equation~\eqref{eq:growth} then gives \eqref{eq:6.18}.
\end{proof}

\section{The two mass constraints when \texorpdfstring{$q\ne0$}{q is nonzero}}\label{sec:noncritical-mass}

Let $u=U+v$ be the projected solution from Theorem~\ref{thm:5.2}. Define
\begin{equation}\label{eq:7.1}
 F_i(a,\lambda,z)=\int_{\R^N}u_i^2-m_i,\qquad i=1,2.
\end{equation}
Both functions are continuous by Lemma~\ref{lem:5.3}. We use the same notation in the critical case.

\begin{proposition}\label{prop:noncritical-mass}
Let $q\ne0$. Uniformly on
$\mathcal R_a\times\overline{\mathcal Z}$,
\begin{equation}\label{eq:7.2}
 \int_{\R^N}u_i^2=k_iM_Wa_i^{-\alpha}\lambda_i^q
 \left(1+\frac{q\overline V_i}{\lambda_i}+O(\lambda_i^{-2})\right),\qquad
\overline V_i=k_i^{-1}\sum_j V_i(b_{i,j}).
\end{equation}
\end{proposition}

\begin{proof}
By rescaling each peak and using the exponential separation of
the peaks within the \(i\)-th component, we obtain
\[
\int_{\mathbb R^N}U_i^2\,dx
=
M_Wa_i^{-\alpha}\sum_{j=1}^{k_i}\omega_{i,j}^{q}
+
O\!\left(a_i^{-\alpha}\lambda_i^q
e^{-c\sqrt{\lambda_i}d_0}\right).
\]
Since \(\xi_{i,j}-b_{i,j}=O(\lambda_i^{-1})\) and
\(\nabla V_i(b_{i,j})=0\), Taylor expansion gives
\[
V_i(\xi_{i,j})=V_i(b_{i,j})+O(\lambda_i^{-2}).
\]
Using \(\omega_{i,j}=\lambda_i+V_i(\xi_{i,j})\), we then have
\[
\omega_{i,j}^{q}
=
\lambda_i^q
\left(1+\frac{qV_i(b_{i,j})}{\lambda_i}
+O(\lambda_i^{-2})\right).
\]
Summing over \(j\) gives the mass expansion for \(U_i\).

The pointwise estimate for \(v_i\) and the frequency growth
condition \eqref{eq:growth} yield
\[
\|v_i\|_{L^2}
\le CA_i\left(\lambda_i^{-2-N/4}
+\tau_\lambda\lambda_*^{-N/4}\right)
\le CA_i\lambda_i^{-2-N/4}.
\]
Since \(\|U_i\|_{L^2}\le CA_i\lambda_i^{-N/4}\) and
\(A_i^2\lambda_i^{-N/2}=a_i^{-\alpha}\lambda_i^q\), we obtain
\[
\begin{aligned}
\left|\int_{\mathbb R^N}(u_i^2-U_i^2)\,dx\right|
\le 2\|U_i\|_{L^2}\|v_i\|_{L^2}+\|v_i\|_{L^2}^2
\le Ca_i^{-\alpha}\lambda_i^{q-2}.
\end{aligned}
\]
Combining this estimate with the mass expansion for \(U_i\)
proves \eqref{eq:7.2}.
\end{proof}

\begin{proposition}\label{prop:noncritical-degree}
For sufficiently large target frequencies there is a point in
$\mathcal Z\times\operatorname{int}\mathcal R_a$ at which
$\mathcal G=0$ and $F_1=F_2=0$.
\end{proposition}

\begin{proof}
Put $I=(1-\eta_0,1+\eta_0)$, with $\eta_0$ as in \eqref{eq:rectangle},
and set $D=\mathcal Z\times I^2$. We regard the center and mass equations
as functions of $(z,t)\in\overline D$ by writing
\[
\lambda_i=t_i\Lambda_i,\qquad \xi_{i,j}=b_{i,j}+\frac{z_{i,j}}{\lambda_i}.
\]
Since $k_iM_Wa_i^{-\alpha}\Lambda_i^q=m_i$, \eqref{eq:7.2} gives
\[
 \frac{F_i}{m_i}=t_i^q-1+o(1)
\]
uniformly on $\overline D$. By \eqref{eq:6.14},
$(\mathcal G,\frac{F_1}{m_1},\frac{F_2}{m_2})$ thus converges uniformly on $\overline D$ to
\[
 (z,t)\longmapsto
 \left((D^2V_i(b_{i,j})(z_{i,j}-z^0_{i,j}))_{i,j},\,
             t_1^q-1,\,t_2^q-1\right).
\]
The only zero of this map in $\overline D$ is $(z^0,1,1)$. The derivative
there has diagonal blocks $D^2V_i(b_{i,j})$, $q$ and $q$, all invertible.
The limiting map is therefore bounded away from zero on $\partial D$ and
has degree $1$ or $-1$. Uniform convergence and homotopy invariance give
the same degree for the center and mass equations.
\end{proof}

\begin{proposition}\label{prop:true-max}
Fix $R_p>1$. Let $u$ be a projected solution for which the center and mass
equations hold. For sufficiently large target frequencies, $u_i$ has
exactly one critical point $x_{i,j}$ in
$B_{\frac{R_p}{\sqrt{\lambda_i}}}(\xi_{i,j})$. It is a strict local maximum, and
\begin{equation}\label{eq:7.16}
 |x_{i,j}-\xi_{i,j}|\le C\lambda_i^{-\frac{5}{2}}.
\end{equation}
Moreover,
\begin{equation}\label{eq:7.17}
 \lim_{R\to\infty}\limsup_{\lambda_*\to\infty}
 \frac{1}{A_i}\sup_{\operatorname{dist}(x,\{\xi_{i,1},\ldots,\xi_{i,k_i}\})
                         \ge \frac{R}{\sqrt{\lambda_i}}} |u_i(x)|=0 .
\end{equation}
\end{proposition}

\begin{proof}
Fix \(i,j\) and set
\[
\xi=\xi_{i,j},\qquad \omega=\omega_{i,j},\qquad
A=a_i^{-1/(p-2)}\omega^{1/(p-2)},\qquad
x(y)=\xi+\frac{y}{\sqrt{\omega}}.
\]
Write \(\widehat u=u_i(x(y))/A\), \(h=\widehat u-W\), and
\(\widehat u_k=u_k(x(y))/A_k\), where \(k=3-i\).
By peak separation and the estimates in
\eqref{eq:peak-decay} and \eqref{eq:5.28}--\eqref{eq:5.30},
for every fixed \(R>0\),
\[
\|h\|_{C^1(B_{2R})}\le C_R\lambda_i^{-2},
\qquad
\|\widehat u_k\|_{C^1(B_{2R})}
\le C_R(1+\lambda^*)^M\tau_\lambda
\]
for some fixed \(M>0\).

Subtracting the equation for \(W\) from the rescaled equation gives
\[
-\Delta h+h
=f(\widehat u)-f(W)-P\widehat u
+\beta\omega^{-1}A^{r-2}A_k^r g(\widehat u,\widehat u_k),
\qquad
P(y)=\frac{V_i(x(y))-V_i(\xi)}{\omega}.
\]
Since \(\nabla V_i(\xi)=O(\lambda_i^{-1})\), Taylor's formula gives
\[
\|P\|_{C^1(B_{2R})}\le C_R\lambda_i^{-2}.
\]
The \(C^1\) estimate for \(h\) gives the same bound for
\(f(\widehat u)-f(W)\). Moreover, \(\widehat u\) is bounded away
from zero on \(B_{2R}\) for large frequencies. The estimate for
\(\widehat u_k\) and \eqref{eq:growth} therefore make the mixed
term \(o(\lambda_i^{-2})\) in \(C^1(B_{2R})\).
The interior Schauder estimate yields
\[
\|h\|_{C^{2,\vartheta}(B_R)}
\le C_{R,\vartheta}\lambda_i^{-2},
\qquad 0<\vartheta<1.
\]

For completeness, the required properties of $W$ follow from the radial energy
\[
E(t)=\frac12W'(t)^2-\frac12W(t)^2+\frac1pW(t)^p,
\qquad E'(t)=-\frac{N-1}{t}W'(t)^2,\qquad E(\infty)=0.
\]
Thus $E\ge0$, and any zero of $W'$ satisfies $W^{p-2}\ge p/2>1$.
The radial equation gives
\[
NW''(0)=W(0)-W(0)^{p-1}<0,\qquad D^2W(0)=W''(0)I.
\]
Hence $W'<0$ near zero. At a first positive zero of $W'$, the same equation would give $W''<0$, a contradiction. Consequently $W'(t)<0$ for all $t>0$.

We now locate the critical point of \(\widehat u\).
Choose \(\delta\in(0,1)\) so that \(D^2W\le-cI\) in \(B_\delta\).
For sufficiently large frequencies, the estimates for \(h\) give
\[
\widehat u(0)>\max_{\partial B_\delta}\widehat u,
\qquad
D^2\widehat u\le-\frac c2 I\quad\text{in }B_\delta.
\]
Thus \(\widehat u\) attains its maximum on \(\overline B_\delta\)
at an interior point \(y_{i,j}\). Strict concavity makes it the
unique critical point in \(B_\delta\) and a strict local maximum.
On the remaining annulus, for sufficiently large frequencies,
\[
|\nabla\widehat u(y)|
\ge \frac12\min_{\delta\le |z|\le 2R_p}|\nabla W(z)|>0,
\qquad \delta\le |y|\le 2R_p.
\]

To prove \eqref{eq:7.16}, use \(\nabla\widehat u(y_{i,j})=0\)
and \(\nabla W(0)=0\) to obtain
\[
\left(\int_0^1D^2W(ty_{i,j})\,dt\right)y_{i,j}
=-\nabla h(y_{i,j}).
\]
Hence \(c|y_{i,j}|\le|\nabla h(y_{i,j})|\le C\lambda_i^{-2}\).
Since \(\omega/\lambda_i\to1\), the rescaled radius
\(R_p\sqrt{\omega/\lambda_i}\) lies between \(\delta\) and \(2R_p\).
Thus \(x_{i,j}=x(y_{i,j})\) is the unique critical point in the
stated ball, and
\[
|x_{i,j}-\xi_{i,j}|=\frac{|y_{i,j}|}{\sqrt{\omega}}
\le C\lambda_i^{-5/2}.
\]

Finally, outside the balls \(B_{R/\sqrt{\lambda_i}}(\xi_{i,j})\),
the pointwise estimates for \(U_i\) and \(v_i\) give
\[
\frac{|u_i(x)|}{A_i}
\le C\bigl(e^{-cR}+\tau_\lambda\bigr).
\]
Letting first \(\lambda_*\to\infty\) and then \(R\to\infty\)
proves \eqref{eq:7.17}.
\end{proof}

\begin{lemma}\label{lem:recenter}
Let $x_{i,j}$ be given by Proposition~\ref{prop:true-max} and put
$\widehat\omega_{i,j}=\lambda_i+V_i(x_{i,j})$. The corresponding recentered
peak $\widehat U_{i,j}=a_i^{-\frac{1}{p-2}}W_{\widehat\omega_{i,j},x_{i,j}}$ satisfies
\begin{equation}\label{eq:7.21}
 |U_{i,j}(x)-\widehat U_{i,j}(x)|
 \le CA_i\lambda_i^{-2}e^{-c\sqrt{\lambda_i}|x-x_{i,j}|}.
\end{equation}
\end{lemma}

\begin{proof}
Write $d=|x_{i,j}-\xi_{i,j}|$. By \eqref{eq:7.16} and (V1), we have
\[
d=O(\lambda_i^{-\frac{5}{2}}),\qquad
|\widehat\omega_{i,j}-\omega_{i,j}|\le\|\nabla V_i\|_\infty d.
\]
Along the parameter segment joining
\((\omega_{i,j},\xi_{i,j})\) and
\((\widehat\omega_{i,j},x_{i,j})\), differentiation with respect to the center and the frequency gives
\[
\begin{aligned}
\left|\nabla_\xi
  \bigl(a_i^{-1/(p-2)}W_{\omega,\xi}(x)\bigr)\right|
&\le C A_i\sqrt{\lambda_i}\,
  e^{-c\sqrt{\lambda_i}|x-x_{i,j}|},\\
\left|\partial_\omega
  \bigl(a_i^{-1/(p-2)}W_{\omega,\xi}(x)\bigr)\right|
&\le C A_i\lambda_i^{-1}\,
  e^{-c\sqrt{\lambda_i}|x-x_{i,j}|}.
\end{aligned}
\]
The linear factor in the scaled distance arising from the frequency
derivative is absorbed into the exponential by decreasing \(c\).

Integrating along the parameter segment, we obtain
\[
\begin{aligned}
|U_{i,j}(x)-\widehat U_{i,j}(x)|
&\le C A_i e^{-c\sqrt{\lambda_i}|x-x_{i,j}|}
 \left(\sqrt{\lambda_i}\,d
 +\lambda_i^{-1}|\widehat\omega_{i,j}-\omega_{i,j}|\right)\\
&\le C A_i\lambda_i^{-2}e^{-c\sqrt{\lambda_i}|x-x_{i,j}|},
\end{aligned}
\]
which proves \eqref{eq:7.21}.
\end{proof}

\begin{proof}[Proof of Theorem~\ref{thm:noncritical}]
At the zero given by Proposition 7.2, (6.5) gives \(c=0\), so the projected solution solves the original system. 
Positivity follows from Section 5.3, and \(F_i=0\) gives the prescribed masses.

Since $\lambda_i\in[(1-\eta_0)\Lambda_i,(1+\eta_0)\Lambda_i]$, we already have $\lambda_i\asymp\Lambda_i$. The mass constraint and \eqref{eq:7.2} give
\[
\begin{aligned}
\Lambda_i
&=\lambda_i\left(1+\frac{q\overline V_i}{\lambda_i}
+O(\lambda_i^{-2})\right)^{1/q}\\
&=\lambda_i+\overline V_i+O(\lambda_i^{-1}).
\end{aligned}
\]
Replacing the last error by $O(\Lambda_i^{-1})$ proves \eqref{eq:frequency-asymptotics}.

Equation~\eqref{eq:6.14} gives
\[
z_{i,j}=z^0_{i,j}+o(1),
\]
so \eqref{eq:7.16} yields \eqref{eq:center-expansion}.
By Lemma~\ref{lem:recenter}, replacing the centers in \eqref{eq:5.28} adds at most the right-hand side of \eqref{eq:7.21} for each peak.
The shifted exponential weights are uniformly equivalent by \eqref{eq:7.16}, and \eqref{eq:growth} gives
\[
\lambda_i^2\tau_\lambda\le e^{-\sigma_0 d_0\sqrt{\lambda_*}}
\]
for large frequencies. Absorbing $\sum_j(U_{i,j}-\widehat U_{i,j})$ into $v_i$ proves
\eqref{eq:solution-profile} and \eqref{eq:correction-bound}. The new correction belongs to $H^2$ by
Theorem~\ref{thm:5.2} and the regularity of the recentered profiles.
Finally, by \eqref{eq:amplitude},
\[
\frac{u_i(x_{i,j})}{A_i}\to W(0),\qquad A_i\to\infty.
\]
Outside any fixed neighborhood of $B_i$, the profile
and the correction tend uniformly to zero by \eqref{eq:growth} and
\eqref{eq:correction-bound}. This verifies Definition~\ref{def:concentrating}.
\end{proof}

\section{The \texorpdfstring{$L^2$}{L2}-critical case}\label{chap:L2-critical}

Suppose $q=0$, so that $(N,p)=(1,6)$ or $(2,4)$. Each peak satisfies
$a_i^\alpha\int U_{i,j}^2=M_W$ for every frequency. To find the next term
we use the scaling derivative
\[
 Z_0=\frac N4 W+\frac12 y\cdot\nabla W,\qquad
 L_0=-\Delta+1-(p-1)W^{p-2}.
\]

\begin{lemma}\label{lem:8.1}
The following identities hold:
\begin{equation}\label{eq:8.2}
 L_0 Z_0=-W,\qquad \int WZ_0=0,\qquad
 \int |y|^2 WZ_0=-\frac12 S_W .
\end{equation}
\end{lemma}

\begin{proof}
Set \(W_s(y)=s^{N/4}W(\sqrt{s}\,y)\), so that
\(\left.\partial_s W_s\right|_{s=1}=Z_0\).
Differentiating
\[
-\Delta W_s+sW_s=W_s^{p-1},\qquad
\int W_s^2=M_W,\qquad
\int |y|^2W_s^2=s^{-1}S_W
\]
at \(s=1\) gives the three identities.
\end{proof}

\noindent\textbf{Remark.}
The identity in Lemma~\ref{lem:8.1} also gives an explicit expression for
the scalar coefficient appearing in the critical mass expansions of
\cite{PellacciPistoiaVairaVerzini,HuangPistoiaTroestlerWang26}.
Let $\Psi\in H^2_{\mathrm{rad}}(\R^N)$ solve
\[
L_0\Psi=|y|^2W\quad\text{in }\R^N.
\]
Such a radial solution exists and is unique by the nondegeneracy of $W$.
Self-adjointness of $L_0$ and \eqref{eq:8.2} yield
\[
\int_{\R^N}W\Psi
=-\int_{\R^N}(L_0Z_0)\Psi
=-\int_{\R^N}Z_0L_0\Psi
=-\int_{\R^N}|y|^2WZ_0
=\frac12S_W>0.
\]
In particular, with the normalization
$\mathfrak m=(2N)^{-1}\int_{\R^N}W\Psi$, one has
$\mathfrak m=S_W/(4N)>0$.
This explains the universal factor $S_W/(2N)$ in
\eqref{eq:critical-coefficient}; its sign is fixed, while the sign of
$\kappa_i$ is determined by $\sum_j\Delta V_i(b_{i,j})$.

\begin{lemma}\label{lem:8.2}
Uniformly on $\mathcal R_a\times\overline{\mathcal Z}$, the coefficients in
\eqref{eq:5.25} satisfy
\begin{equation}\label{eq:8.3}
 \max_{j,l}|c_{i,j,l}|\le C\lambda_i^{-\frac{1}{2}}.
\end{equation}
\end{lemma}

\begin{proof}
By \eqref{eq:6.6} and $|\nabla V_i(\xi_{i,j})|=O(\lambda_i^{-1})$,
$|G_{i,j,l}|\le Cs_i\lambda_i^{-1}$. Now use
\eqref{eq:6.2} and \eqref{eq:6.4}.
\end{proof}

Choose \(\rho>0\) with \(2\rho<d_0\) so that, for large frequencies,
\(B_{2\rho}(\xi_{i,j})\) lies in the corresponding neighborhood in (V2).
Let
\[
\chi\in C_c^\infty(B_{2\rho}(0)),\qquad
0\le\chi\le1,\qquad
\chi\equiv1 \text{in } B_\rho(0),\qquad
\chi_{i,j}(x)=\chi(x-\xi_{i,j}).
\]
Fix \(i,j\) and write
\[
\xi=\xi_{i,j},\qquad \omega=\omega_{i,j},\qquad
A=a_i^{-1/(p-2)}\omega^{1/(p-2)},\qquad
x(y)=\xi+\frac{y}{\sqrt{\omega}}.
\]
Define
\[
\chi_\omega(y)=\chi\left(\frac{y}{\sqrt{\omega}}\right),\qquad
\widehat u(y)=\frac{u_i(x(y))}{A},\qquad
\widetilde u=\chi_\omega\widehat u,\qquad h=\widetilde u-W.
\]
Since \(a_i^\alpha A^2\omega^{-N/2}=\omega^q=1\), changing variables gives
\begin{equation}\label{eq:8.5}
a_i^\alpha\int\chi_{i,j}^2u_i^2=\int\widetilde u^2.
\end{equation}

\begin{lemma}\label{lem:8.3}
Uniformly on $\mathcal R_a\times\overline{\mathcal Z}$,
\begin{equation}\label{eq:8.6}
 \|h\|_\infty+\|h\|_2\le C\lambda_i^{-2},\qquad
 a_i^\alpha\int\Bigl(1-\sum_j\chi_{i,j}^2\Bigr)u_i^2
 =o(\lambda_i^{-2}).
\end{equation}
With $P(y)=\frac{V_i(x(y))-V_i(\xi)}{\omega}$, one has
\begin{equation}\label{eq:8.7}
 -\Delta\widetilde u+(1+P)\widetilde u=f(\widetilde u)+E,\qquad
 \int|Z_0E|=o(\lambda_i^{-2}).
\end{equation}
\end{lemma}

\begin{proof}
Since \(u_i=\sum_s U_{i,s}+v_i\) and \(U_{i,j}(x(y))=AW(y)\), we have
\[
h=(\chi_\omega-1)W+
\frac{\chi_\omega}{A}
\left(\sum_{s\ne j}U_{i,s}+v_i\right)(x(y)).
\]
The exponential decay of \(W\) and peak separation give
\[
\|(\chi_\omega-1)W\|_\infty+\|(\chi_\omega-1)W\|_2
\le Ce^{-c\sqrt{\omega}}
\]
and
\[
\left\|\chi_\omega\sum_{s\ne j}\frac{U_{i,s}(x(\cdot))}{A}\right\|_\infty
+\left\|\chi_\omega\sum_{s\ne j}\frac{U_{i,s}(x(\cdot))}{A}\right\|_2
\le Ce^{-c\sqrt{\omega}}.
\]
By \eqref{eq:5.28}, \(A_i/A=1+O(\lambda_i^{-1})\), and the change
of variables \(dy=\omega^{N/2}\,dx\),
\[
\left\|\frac{\chi_\omega v_i(x(\cdot))}{A}\right\|_\infty
+\left\|\frac{\chi_\omega v_i(x(\cdot))}{A}\right\|_2
\le C\left[\lambda_i^{-2}
+\tau_\lambda\left(\frac{\lambda_i}{\lambda_*}\right)^{N/4}\right]
\le C\lambda_i^{-2},
\]
where the last inequality follows from \eqref{eq:growth}.
This proves the first estimate in \eqref{eq:8.6}.

Let \(\Omega_{\mathrm{out}}\) be the support of
\(1-\sum_j\chi_{i,j}^2\). Since the cutoffs have disjoint supports,
\(\Omega_{\mathrm{out}}\) lies outside the balls
\(B_\rho(\xi_{i,j})\).
The exponential tails of the principal peaks and \eqref{eq:5.28} give
\[
a_i^\alpha\int_{\Omega_{\mathrm{out}}}u_i^2
\le Ce^{-c\sqrt{\lambda_i}}
+C\tau_\lambda^2\left(\frac{\lambda_i}{\lambda_*}\right)^{N/2}
=o(\lambda_i^{-2})
\]
by \eqref{eq:growth}. Since \(0\le1-\sum_j\chi_{i,j}^2\le1\),
this proves the second estimate in \eqref{eq:8.6}.

Let \(k=3-i\) and set \(\widehat u_k(y)=u_k(x(y))/A_k\).
Multiplying the rescaled projected equation by \(\chi_\omega\)
gives the equation in \eqref{eq:8.7}, with
\[
\begin{aligned}
E
&=\underbrace{
  -2\nabla\chi_\omega\cdot\nabla\widehat u
  -(\Delta\chi_\omega)\widehat u
  +\chi_\omega f(\widehat u)-f(\chi_\omega\widehat u)
}_{E_1}\\
&\quad+\underbrace{
  \chi_\omega\frac{\beta A^{r-2}A_k^r}{\omega}
  g(\widehat u,\widehat u_k)
}_{E_2}
+\underbrace{
  \frac{\chi_\omega}{\omega A}
  \sum_{s,l}c_{i,s,l}(D_iY_{i,s,l})(x(y))
}_{E_3}.
\end{aligned}
\]
By the support of \(E_1\), peak separation, and
\eqref{eq:growth}, \eqref{eq:4.12}, and \eqref{eq:8.3},
\[
\int |Z_0E_1|\le Ce^{-c\sqrt{\omega}},\qquad
\int |Z_0E_2|=o(\lambda_i^{-2}),\qquad
\int |Z_0E_3|\le C\lambda_i^{-5/2}.
\]
Thus \(\int |Z_0E|=o(\lambda_i^{-2})\).
\end{proof}

\begin{proposition}\label{prop:8.4}
Uniformly on $\mathcal R_a\times\overline{\mathcal Z}$,
\begin{equation}\label{eq:8.10}
 a_i^\alpha\int_{\R^N}u_i^2
 =k_iM_W-\frac{\kappa_i}{\lambda_i^2}+o(\lambda_i^{-2}),\qquad
 \kappa_i=\frac{S_W}{2N}\sum_{j=1}^{k_i}\Delta V_i(b_{i,j}).
\end{equation}
\end{proposition}

\begin{proof}

Fix a peak. Since \(\widetilde u=W+h\), \eqref{eq:8.6} gives
\[
\int\widetilde u^2=M_W+2\int Wh+O(\lambda_i^{-4}).
\]
Subtracting the equation for \(W\) from \eqref{eq:8.7} yields
\[
L_0h=-P(W+h)+Q(h)+E,
\qquad
Q(h)=f(W+h)-f(W)-f'(W)h.
\]
By the self-adjointness of \(L_0\) and \(L_0Z_0=-W\) in
\eqref{eq:8.2},
\[
\int Wh
=\int Z_0PW+\int Z_0Ph-\int Z_0Q(h)-\int Z_0E.
\]
Since \(\|P\|_\infty=O(\lambda_i^{-1})\),
\(\|h\|_\infty=O(\lambda_i^{-2})\), and
\(\|Q(h)\|_\infty=O(\|h\|_\infty^2)\),
the estimates in \eqref{eq:8.6}--\eqref{eq:8.7} give
\[
\int Wh=\int Z_0PW+o(\lambda_i^{-2}).
\]

For \(|y|\le\rho\sqrt{\omega}\), Taylor's formula gives
\[
P(y)=\frac{\nabla V_i(\xi)\cdot y}{\omega^{3/2}}
+\frac{y^TD^2V_i(\xi)y}{2\omega^2}
+O(\omega^{-5/2}|y|^3).
\]
The exponential decay of \(WZ_0\) makes the integral outside this ball
exponentially small and bounds the integrated Taylor remainder by
\(O(\omega^{-5/2})\). The linear term integrates to zero by oddness.
For the quadratic term, radial symmetry and \eqref{eq:8.2} give
\[
\int y_l y_m WZ_0
=\frac{\delta_{lm}}{N}\int |y|^2WZ_0
=-\frac{\delta_{lm}S_W}{2N}.
\]
Thus
\[
\begin{aligned}
\int Z_0PW
&=-\frac{S_W}{4N\omega^2}\Delta V_i(\xi)
  +o(\lambda_i^{-2})\\
&=-\frac{S_W}{4N\lambda_i^2}\Delta V_i(b_{i,j})
  +o(\lambda_i^{-2}).
\end{aligned}
\]
The last equality follows from \(\omega=\lambda_i+O(1)\) and
\(\xi=b_{i,j}+O(\lambda_i^{-1})\), together with (V2).

Summing the local mass expansions using \eqref{eq:8.5} and the
exterior mass estimate in \eqref{eq:8.6} proves \eqref{eq:8.10}.
\end{proof}

\begin{proof}[Proof of Theorem~\ref{thm:critical}]

Let \(a=(a_1,a_2)\) and \(t_i=\lambda_i/\Lambda_i(a_i)\). Use the fixed
domain from the proof of Proposition~\ref{prop:noncritical-degree},
\[
D=\mathcal Z\times(1-\eta_0,1+\eta_0)^2,
\qquad
\lambda(t)=\bigl(t_1\Lambda_1(a_1),t_2\Lambda_2(a_2)\bigr).
\]
By \eqref{eq:critical-direction} and \eqref{eq:critical-target},
\(\delta_i(a_i)=\kappa_i/\Lambda_i(a_i)^2\ne0\). The mass expansion
\eqref{eq:8.10} gives
\[
\frac{a_i^\alpha F_i(a,\lambda(t),z)}{\delta_i(a_i)}
=1-t_i^{-2}+o(1),\qquad i=1,2.
\]
The remainder is uniform on \(\overline D\), since
\(
\frac{\lambda_i^{-2}}{|\delta_i(a_i)|}
=\frac{t_i^{-2}}{|\kappa_i|}
\)
is bounded there. Define the joint map on \(\overline D\) by
\[
\mathcal H_a^{\mathrm c}(z,t_1,t_2)
:=
\left(
\mathcal G(a,\lambda(t),z),\,
\frac{a_1^\alpha F_1(a,\lambda(t),z)}{\delta_1(a_1)},\,
\frac{a_2^\alpha F_2(a,\lambda(t),z)}{\delta_2(a_2)}
\right).
\]
By \eqref{eq:6.14} and the mass expansion above,
\(\mathcal H_a^{\mathrm c}\) converges uniformly on \(\overline D\) to
\[
\mathcal H_0^{\mathrm c}(z,t_1,t_2)
:=
\left(
\bigl(D^2V_i(b_{i,j})(z_{i,j}-z_{i,j}^0)\bigr)_{i,j},\,
1-t_1^{-2},\,1-t_2^{-2}
\right).
\]
The unique zero of \(\mathcal H_0^{\mathrm c}\) in \(\overline D\) is
\((z^0,1,1)\in D\). Its derivative there has invertible diagonal blocks
\(D^2V_i(b_{i,j})\), \(2\), and \(2\), so the zero is nondegenerate.
The Brouwer degree argument in the proof of
Proposition~\ref{prop:noncritical-degree} gives
\[
\deg(\mathcal H_a^{\mathrm c},D,0)
=\deg(\mathcal H_0^{\mathrm c},D,0)\ne0
\]
for sufficiently large target frequencies. Hence
\(\mathcal H_a^{\mathrm c}\) has a zero in \(D\), at which
\[
\mathcal G(a,\lambda(t),z)=0,\qquad
F_1(a,\lambda(t),z)=F_2(a,\lambda(t),z)=0.
\]
At this zero, the normalized mass equations give \(1-t_i^{-2}=o(1)\).
Thus \(t_i\to1\), which proves \eqref{eq:critical-frequency}.

By \eqref{eq:6.5}, \(\mathcal G=0\) makes the projection coefficients
vanish, so the solution satisfies the original system. Positivity follows
from Section~\ref{sec:positivity}, and \(F_1=F_2=0\) gives the prescribed
masses. The peak location expansion, the recentered profile decomposition
and correction estimate, and the segregated multi-peak concentration
follow as in the proof of Theorem~\ref{thm:noncritical}.
\end{proof}

\appendix

\section{Cutoff errors in the Pohozaev identities}\label{app:nehari-pohozaev}
For the global identity, take \(\chi_R=\chi(\frac{x}{R})\), with \(\chi\in C_c^\infty(B_2)\) equal to one on \(B_1\), and test the two equations by \(\chi_R x\cdot\nabla u_i\). The sum of the mixed terms is
\[
\frac{\beta}{r}\int \chi_R x\cdot\nabla(u_1^r u_2^r)
=-\frac{\beta}{r}\int (N\chi_R+x\cdot\nabla\chi_R)u_1^r u_2^r.
\]
On the support of \(\nabla\chi_R\) we have \(|x||\nabla\chi_R|\le C\).
For fixed \(n\), all terms involving its derivative are therefore bounded by a constant times
\[
\int_{B_{2R}\setminus B_R}\sum_{i=1}^2
\bigl(|\nabla u_i|^2+u_i^2+u_i^p\bigr)\,dx .
\]
We used \(2u_1^ru_2^r\le u_1^p+u_2^p\). The annular integral tends
to zero because \(u_i\in H^1\cap L^p\). Also
\((x\cdot\nabla V_i)u_i^2\in L^1\) by (V1). We can therefore let
\(R\to\infty\) to obtain \eqref{eq:global-pohozaev}.

For the localized identity, take
\(X_\ell=\chi_{i,\ell}^2(x-b_{i,\ell})\) and suppress the index \(n\).
The terms in which a derivative hits a cutoff are supported on
\(T_i:=\bigcup_\ell\operatorname{supp}\nabla\chi_{i,\ell}\).
If \(\mu_{i,n},a_{i,n}\) are bounded, the sum of their absolute values is bounded by
\[
C\int_{T_i} (|\nabla u_i|^2+u_i^2+a_i u_i^p)\,dx=o(1)
\]
by Lemma~\ref{lem:local-energy}. The mixed term is kept as
\(\beta\sum_\ell\int u_i^{r-1}u_{3-i}^r X_\ell\cdot\nabla u_i=(\frac{\beta}{r})J_i^{\mathrm{loc}}\).
The remaining terms give
\[
\tfrac{N-2}{2} K_i^{\rm loc}+\tfrac N2 B_i^{\rm loc}
+\tfrac12\mathcal V_i^{\rm loc}
=\tfrac Np a_i P_i^{\rm loc}-\tfrac{\beta}{r}J_i^{\rm loc}+o(1).
\]
Eliminating \(a_iP_i^{\rm loc}\) by \eqref{eq:local-nehari} gives
\eqref{eq:local-pohozaev}.

\section{The bounded kernel of the scalar linearization}\label{app:scalar-identities}

\begin{lemma}\label{lem:appendix-bounded-kernel}
Every bounded distributional solution of
\[
 L_0\phi:= -\Delta\phi+\phi-(p-1)W^{p-2}\phi=0 \quad\text{in }\R^N
\]
belongs to $\operatorname{span}\{\partial_1W,\ldots,\partial_NW\}$.
\end{lemma}
\begin{proof}
Since $W$ decays exponentially and $\phi$ is bounded,
$F=(p-1)W^{p-2}\phi\in L^2(\R^N)$. Taking the Fourier transform of
$(-\Delta+1)\phi=F$ in tempered distributions gives
$(1+|\zeta|^2)\widehat\phi=\widehat F$, hence $\phi\in H^2(\R^N)$.
The assertion follows from the $L^2$ kernel characterization in
\cite[Lemma~2.1]{ChangGustafsonNakanishiTsai}.
\end{proof}

\noindent\begin{minipage}{\textwidth}

  \medskip
  \noindent\textbf{Acknowledgments.}
  Qing Guo was supported by the National Natural Science Foundation of China (Grant No.~12271539). The authors acknowledge the assistance of ChatGPT (OpenAI) with the language polishing and take full responsibility for the
  mathematical content and the exposition
  of this paper.
  \end{minipage}


\begin{thebibliography}{99}

\bibitem{AmbrosettiBadialeCingolani}
A. Ambrosetti, M. Badiale and S. Cingolani,
Semiclassical states of nonlinear Schr\"odinger equations,
\emph{Arch. Ration. Mech. Anal.} \textbf{140} (1997), 285--300.

\bibitem{BaoCai}
W. Bao and Y. Cai,
Mathematical theory and numerical methods for Bose--Einstein condensation,
\emph{Kinet. Relat. Models} \textbf{6} (2013), 1--135.

\bibitem{BartschJeanjean18}
T. Bartsch and L. Jeanjean,
Normalized solutions for nonlinear Schr\"odinger systems,
\emph{Proc. Roy. Soc. Edinburgh Sect. A} \textbf{148} (2018), 225--242.
\newblock doi: \href{https://doi.org/10.1017/S0308210517000087}{10.1017/S0308210517000087}.

\bibitem{BartschSoave17}
T. Bartsch and N. Soave,
A natural constraint approach to normalized solutions of nonlinear
Schr\"odinger equations and systems,
\emph{J. Funct. Anal.} \textbf{272} (2017), 4998--5037.
\newblock doi: \href{https://doi.org/10.1016/j.jfa.2017.01.025}{10.1016/j.jfa.2017.01.025}.

\bibitem{BartschSoaveCorrection}
T. Bartsch and N. Soave,
Correction to: \textquotedblleft A natural constraint approach to normalized solutions of
nonlinear Schr\"odinger equations and systems\textquotedblright\ [J. Funct. Anal. 272 (12) (2017) 4998--5037],
\emph{J. Funct. Anal.} \textbf{275} (2018), 516--521.
\newblock doi: \href{https://doi.org/10.1016/j.jfa.2018.02.007}{10.1016/j.jfa.2018.02.007}.

\bibitem{BartschSoave19}
T. Bartsch and N. Soave,
Multiple normalized solutions for a competing system of Schr\"odinger equations,
\emph{Calc. Var. Partial Differential Equations} \textbf{58} (2019), Article 22.
\newblock doi: \href{https://doi.org/10.1007/s00526-018-1476-x}{10.1007/s00526-018-1476-x}.

\bibitem{BartschZhongZou21}
T. Bartsch, X. Zhong and W. Zou,
Normalized solutions for a coupled Schr\"odinger system,
\emph{Math. Ann.} \textbf{380} (2021), 1713--1740.
\newblock doi: \href{https://doi.org/10.1007/s00208-020-02000-w}{10.1007/s00208-020-02000-w}.

\bibitem{CaoNoussairYan}
D. Cao, E. Noussair and S. Yan,
Solutions with multiple peaks for nonlinear elliptic equations,
\emph{Proc. Roy. Soc. Edinburgh Sect. A} \textbf{129} (1999), 235--264.

\bibitem{CaoPengYan}
D. Cao, S. Peng and S. Yan,
\emph{Singularly Perturbed Methods for Nonlinear Elliptic Problems},
Cambridge University Press, Cambridge, 2021.

\bibitem{ChangGustafsonNakanishiTsai}
S.-M. Chang, S. Gustafson, K. Nakanishi and T.-P. Tsai,
Spectra of linearized operators for NLS solitary waves,
\emph{SIAM J. Math. Anal.} \textbf{39} (2008), 1070--1111.
\newblock doi: \href{https://doi.org/10.1137/050648389}{10.1137/050648389}.

\bibitem{ChenZou12}
Z. Chen and W. Zou,
Positive least energy solutions and phase separation for coupled
Schr\"odinger equations with critical exponent,
\emph{Arch. Ration. Mech. Anal.} \textbf{205} (2012), 515--551.

\bibitem{ChenZou15}
Z. Chen and W. Zou,
Positive least energy solutions and phase separation for coupled
Schr\"odinger equations with critical exponent: higher dimensional case,
\emph{Calc. Var. Partial Differential Equations} \textbf{52} (2015),
423--467.

\bibitem{ClappP18}
M. Clapp and A. Pistoia,
Existence and phase separation of entire solutions to a pure critical
competitive elliptic system,
\emph{Calc. Var. Partial Differential Equations} \textbf{57} (2018),
Article 23.

\bibitem{GidasSpruck}
B. Gidas and J. Spruck,
Global and local behavior of positive solutions of nonlinear elliptic
equations,
\emph{Comm. Pure Appl. Math.} \textbf{34} (1981), 525--598.

\bibitem{GuoHua25}
Q. Guo and Q. Hua,
Segregated solutions for a nonlinear Schr\"odinger system involving mass
supercritical exponents,
\emph{Nonlinearity} \textbf{38} (2025), 035009.
\newblock doi: \href{https://doi.org/10.1088/1361-6544/adb196}{10.1088/1361-6544/adb196}.

\bibitem{GuoYang23}
Q. Guo and J. Yang,
Excited states for two-component Bose--Einstein condensates in dimension two,
\emph{J. Differential Equations} \textbf{343} (2023), 659--686.

\bibitem{GuoZhangCV}
Q. Guo and C. Zhang,
Segregated solutions for nonlinear Schr\"odinger systems with sublinear coupling terms,
\emph{Calc. Var. Partial Differential Equations} \textbf{65} (2026), Article 70.
\newblock doi: \href{https://doi.org/10.1007/s00526-025-03239-4}{10.1007/s00526-025-03239-4}.

\bibitem{GuoSeiringer}
Y. Guo and R. Seiringer,
On the mass concentration for Bose--Einstein condensates with attractive
interactions,
\emph{Lett. Math. Phys.} \textbf{104} (2014), 141--156.

\bibitem{HuangPistoiaTroestlerWang26}
X. Huang, A. Pistoia, C. Troestler and C. Wang,
Normalized vector solutions of nonlinear Schr\"odinger systems,
\emph{Discrete Contin. Dyn. Syst.} \textbf{53} (2026), 188--209.
\newblock doi: \href{https://doi.org/10.3934/dcds.2025162}{10.3934/dcds.2025162}.

\bibitem{HuLiLongWang26}
W. Hu, B. Li, W. Long and C. Wang,
Existence and uniqueness of normalized multi-peak solutions for coupled
nonlinear Schr\"odinger systems,
preprint (2026), \href{https://arxiv.org/abs/2604.27455}{arXiv:2604.27455}.

\bibitem{Kwong}
M. K. Kwong,
Uniqueness of positive solutions of $\Delta u-u+u^p=0$ in $\R^N$,
\emph{Arch. Ration. Mech. Anal.} \textbf{105} (1989), 243--266.

\bibitem{LuoPengWeiYan}
P. Luo, S. Peng, J. Wei and S. Yan,
Excited states of Bose--Einstein condensates with degenerate attractive
interactions,
\emph{Calc. Var. Partial Differential Equations} \textbf{60} (2021),
Article 155.

\bibitem{MederskiSzulkin}
J. Mederski and A. Szulkin,
Multiple normalized solutions to a system of nonlinear Schr\"odinger equations,
\emph{J. Geom. Anal.} \textbf{35} (2025), Article 312.

\bibitem{NorisTavaresVerzini}
B. Noris, H. Tavares and G. Verzini,
Stable solitary waves with prescribed $L^2$-mass for the cubic Schr\"odinger
system with trapping potentials,
\emph{Discrete Contin. Dyn. Syst.} \textbf{35} (2015), 6085--6112.

\bibitem{PellacciPistoiaVairaVerzini}
B. Pellacci, A. Pistoia, G. Vaira and G. Verzini,
Normalized concentrating solutions to nonlinear elliptic problems,
\emph{J. Differential Equations} \textbf{275} (2021), 882--919.
\newblock doi: \href{https://doi.org/10.1016/j.jde.2020.11.003}{10.1016/j.jde.2020.11.003}.

\bibitem{PengWang13}
S. Peng and Z.-Q. Wang,
Segregated and synchronized vector solutions for nonlinear Schr\"odinger
systems,
\emph{Arch. Ration. Mech. Anal.} \textbf{208} (2013), 305--339.

\bibitem{TangSunSun26}
W.-Y. Tang, H.-R. Sun and G.-P. Sun,
Multiple normalized solutions for a class of non-autonomous
Schr\"odinger system in $\R^2$,
\emph{J. Math. Anal. Appl.} \textbf{556} (2026), 130255.
\newblock doi: \href{https://doi.org/10.1016/j.jmaa.2025.130255}{10.1016/j.jmaa.2025.130255}.

\bibitem{TavaresTerracini12}
H. Tavares and S. Terracini,
Sign-changing solutions of competition--diffusion elliptic systems and
optimal partition problems,
\emph{Ann. Inst. H. Poincar\'e Anal. Non Lin\'eaire}
\textbf{29} (2012), 279--300.

\bibitem{Weinstein}
M. I. Weinstein,
Modulational stability of ground states of nonlinear Schr\"odinger equations,
\emph{SIAM J. Math. Anal.} \textbf{16} (1985), 472--491.

\end{thebibliography}
\end{document}